\documentclass[10pt]{amsart}
\usepackage{graphicx}
\usepackage{amsmath}
\usepackage{amsthm}
\usepackage{amsfonts}
\usepackage{amssymb}
\usepackage{cite}

\usepackage{pstricks}

\newtheorem{theorem}{Theorem}[section]
\newtheorem{lemma}[theorem]{Lemma}
\newtheorem{proposition}[theorem]{Proposition}
\newtheorem{corollary}[theorem]{Corollary}

\theoremstyle{definition}
\newtheorem{definition}[theorem]{Definition}

\theoremstyle{remark}
\newtheorem{remark}[theorem]{Remark}

\numberwithin{equation}{section}

\newcommand{\lp}{\left(}
\newcommand{\rp}{\right)}

\newcommand{\zero}{\overrightarrow{0}}
\newcommand{\nn}{\overrightarrow{n}}
\newcommand{\mm}{\overrightarrow{m}}
\newcommand{\ee}{\overrightarrow{e}}
\renewcommand{\aa}{\overrightarrow{a}}
\newcommand{\n}[1]{\overrightarrow{n_{#1}}}
\newcommand{\m}[1]{\overrightarrow{m_{#1}}}
\renewcommand{\v}[2]{\overrightarrow{V_{#1}V_{#2}}}
\newcommand{\p}[2]{\overrightarrow{P_{#1}P_{#2}}}
\newcommand{\intr}[2]{\overline{#1,#2}}
\newcommand{\interi}{\operatorname{int}}

\renewcommand{\|}{\,\rule[-3.5pt]{.5pt}{12pt}\,}
\newcommand{\notOneSide}{
{\hspace*{3pt}\slash\hspace*{-3pt}{\rule[-3.5pt]{.5pt}{12pt}}\ }
}

\newcommand{\al}{\alpha}
\newcommand{\la}{\lambda}
\newcommand{\si}{\sigma}
\newcommand{\de}{\delta}

\renewcommand{\th}{\theta}

\renewcommand{\le}{\leqslant}
\renewcommand{\ge}{\geqslant}

\newcommand{\R}{\mathbb{R}}
\newcommand{\Z}{\mathbb{Z}}

\renewcommand{\P}{\mathcal{P}}
\newcommand{\Q}{\mathcal{Q}}

\newcommand{\conv}{\operatorname{conv}}
\newcommand{\ext}{\operatorname{ext}}
\newcommand{\edg}{\operatorname{edg}}
\newcommand{\ri}{\operatorname{ri}\,}
\newcommand{\card}{\operatorname{card}}

\begin{document}

\title[Convexity of Sub-polygons of Convex Polygons]{Convexity of Sub-polygons of Convex Polygons}

\author{Iosif Pinelis}
\address{Department of Mathematical Sciences,
Michigan Technological University,
Houghton, MI 49931}
\email{ipinelis@mtu.edu}



\date{\today. File: { ipinelis/polygon/convex-poly/sub-polygons/\jobname.tex
}
}



\begin{abstract}
A convex polygon is defined as a sequence $(V_0,\dots,V_{n-1})$ of points on a plane such that 
the union of the edges 
$[V_0,V_1],\dots,[V_{n-2},V_{n-1}],[V_{n-1},V_0]$
coincides with the boundary of the convex hull of the set of vertices $\{V_0,\dots,V_{n-1}\}$.
It is proved that all sub-polygons of any convex polygon with distinct vertices are convex. 
It is also proved that, if all sub-$(n-1)$-gons of an $n$-gon with $n\ge5$ are convex, then the $n$-gon is convex. 
Other related results are given.
\end{abstract}

\maketitle



\tableofcontents

\setcounter{section}{-1}

\section{Introduction}



Everyone knows a convex polygon when one sees it. However, to deal with the notion of polygon convexity mathematically or computationally, it must be adequately described.  
A convex polygon can be defined, as e.g. in \cite[page 5]{yaglom}, as 
a succession of connected line segments which constitute the boundary of a convex set.
However, in computational geometry it seems more convenient to
consider a polygon as a sequence of its vertices, say $(V_0,\dots,V_{n-1})$, with the edges being the segments $[V_0,V_1],\dots,[V_{n-2},V_{n-1}],[V_{n-1},V_0]$. 
Then one can say that a polygon is convex if the union of its edges coincides with the boundary of the convex hull of the set of vertices $\{V_0,\dots,V_{n-1}\}$.
 
Now let us look at the following picture.
It suggests that, if any one of the vertices of a convex polygon is removed, then the so reduced polygon inherits the 

\vspace*{-5.1cm}

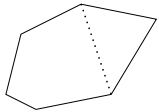
\parbox{2cm}
{
\pspicture(.4,-5.5)(11,7)
\psset{unit=2mm}
\psset{linewidth=.1mm}
\pspolygon(3,0)(8,1)(11,6)(6,7)(2,5)(1,1)
\psline[linewidth=.25mm,linestyle=dotted,dotsep=.3](8,1)(6,7)
\endpspicture
}
\parbox{9.9cm}
{convexity property. Clearly, such a reduction in the number of vertices, $n$, should be helpful, as it could be used to prove various characterizations of convex polygons by induction in $n$.}

\vspace{-5.1cm}

\noindent In particular, the hereditariness of polygon convexity under vertex elimination could be used to establish incremental tests for polygon convexity; cf. \cite{test}.
One finds the following statement in \cite[page 233]{moret}:
\begin{quote}
{\bf Theorem 4.3}\quad Let the sequence of vertices, $p_1,p_2,\dots,p_n,p_{n+1}=p_1$, define an arbitrary polygon $P$ and let $P_i$ be the polygon defined by the sequence of vertices $p_1,p_2,\dots,p_i,p_1$. Then $P$ is convex if and only if, for each $i$,\quad $i=3,4,\dots,n$, polygon $P_i$ is itself convex. 
\end{quote}
It is also said in \cite{moret} that an
incremental test for polygon convexity can be based on the quoted theorem. No proof or reference to a proof of this theorem was given there. Moreover, the ``if" part of the theorem is trivial: if all polygons $P_3,\dots,P_n$ are convex, then polygon $P=P_n$ is trivially convex. Thus, 
it would be impossible to base an incremental test on such a theorem by itself. 

One might suppose that there was a typo in the quoted statement of Theorem~4.3, and there was meant to be $i=3,4,\dots,n-1$ in place of $i=3,4,\dots,n$ (or, equivalently, $p_1,p_2,\dots,p_n,p_{n+1},p_1$ in place of $p_1,p_2,\dots,p_n,p_{n+1}=p_1$). 
But then
the theorem could not be true. 
Indeed, note that all $n$-gons with $n\le3$ are convex. Hence, if the ``if" part of quoted Theorem 4.3 were true with $i=3,4,\dots,n-1$ in place of $i=3,4,\dots,n$, then it would immediately follow by induction in $n$ that all polygons whatsoever are convex! 

As for ``only if" part of Theorem 4.3, whether it is true or not for an arbitrary polygon (given by an arbitrary sequence of vertices $p_1,p_2,\dots,p_n,p_{n+1}=p_1$) 
depends on what definition of polygon convexity is assumed. 
For example,
let $n=8$, $p_1=(0,0)$, $p_2=(1,0)$, $p_3=(1,1)$, $p_4=(0,1)$, $p_5=(0,0)$, $p_6=(1,0)$, $p_7=(1,1)$, $p_8=(0,1)$, $p_{8+1}=p_1=(0,0)$. This polygon ``traces out" the edges of the unit square twice, counter-clockwise. Then this polygon is convex, according to the definition given in the first paragraph of our paper. Yet, the reduced polygon $P_7$, given by $p_1,p_2,\dots,p_7,p_1$ (with vertex $p_8$ removed), is not convex. 

(Other ``obvious'' but incorrect or incomplete statements or proofs were discussed in \cite{jgeom,angles}.)

However, the main result of our paper (Theorem \ref{th:conv-reduction}) states that 
if $\P=(V_0,\dots,V_{n-1})$ is
a convex polygon whose vertices are all {\em distinct}, then the reduced polygon $\P^{(i)}:=(V_0,\dots,V_{i-1},V_{i+1},\dots,V_{n-1})$ (with vertex $V_i$ and hence edges $[V_{i-1},V_i]$ and $[V_i,V_{i+1}]$ removed) is also convex, for each $i$.  

In addition to such downward hereditariness of polygon convexity, we show (Theorem \ref{th:conv-upwards}) that the polygon convexity property is hereditary upwards as well. Namely, if a polygon $\P=(V_0,\dots,V_{n-1})$ with $n\ge5$ vertices (which do not have to be distinct here) is such that all the reduced polygons $\P^{(i)}$ are convex, then
$\P$ is also convex.

It should be clear that the downward hereditariness of polygon convexity given by Theorem \ref{th:conv-reduction} can be used to prove by induction in $n$ that a given condition (say C) is necessary for the polygon convexity, provided that condition C is hereditary upwards. Indeed, suppose that an $n$-gon $\P=(V_0,\dots,V_{n-1})$ is convex, while a condition C is hereditary upwards and satisfied by all convex $k$-gons with $k\le n-1$. Then, by the downward hereditariness of polygon convexity, all the reduced polygons $\P^{(i)}$ are convex. Hence, by induction, all the $\P^{(i)}$'s satisfy condition C. Then the upward hereditariness of condition $C$ will imply that polygon $\P$ also satisfies condition C. 

Similarly, the upward hereditariness of polygon convexity given by Theorem~\ref{th:conv-upwards} can be used to prove by induction that a given condition C is sufficient for the polygon convexity, provided that C is hereditary downwards. 

Thus, taken together, Theorems \ref{th:conv-reduction} and \ref{th:conv-upwards} can be used to obtain conditions necessary and sufficient for polygon convexity. In particular, the test for polygon convexity given by Corollary \ref{cor:test} is immediate from Theorems \ref{th:conv-reduction} and \ref{th:conv-upwards}. Namely, a polygon $\P=(V_0,\dots,V_{n-1})$ with $n\ge5$ distinct vertices is convex if and only if all the reduced polygons $\P^{(i)}$ are convex. 

Such a test should be helpful in theoretical considerations.
However, the test based on a straightforward application of Corollary \ref{cor:test} would be extremely wasteful computationally. 

Indeed, suppose that for every $n\ge5$ one tests the convexity of polygon of polygon $\P=(V_0,\dots,V_{n-1})$ by testing the convexity of all the $n$ reduced polygons $\P^{(i)}$. Then one has $A_n=nA_{n-1}$, where $A_n$ stands for the number of operations needed to test the convexity of an $n$-gon. Hence, $A_n=c\,n!$, for some natural constant $c$. 

A test of a much smaller computational complexity, $O(n^4)$, can be obviously based on Corollary \ref{cor:4-gons}. 

An $O(n^2)$ polygon convexity test---only for simple polygons---could be based on unimodality properties stated as \cite[Theorems 1 and 2]{tous}. (However, the proofs in \cite{tous} are rather heuristic.)
 
Elsewhere \cite{test} we develop and present 
an $O(n)$ polygon convexity test, which is moreover minimal in a certain sense. 

What has been said calls for a rigorous approach---to definitions, statements, and proofs. This is what we aim to provide in this paper. Even though the fact stated in the ``downward'' theorem (Theorem \ref{th:conv-reduction}) seems quite intuitive, the rigorous proof of the theorem (which is the simplest and only proof known to this author) is rather complicated. It is based on a series of 15 lemmas, with numerous logical connections between them. 
On reflection, perhaps it should not be surprising that the proof of Theorem \ref{th:conv-reduction} is complicated. One should realize that the very definition of a convex polygon is rather complex, and it is certainly significantly more complex than the usual definition of a convex set (as a set containing the entire segment whenever the set contains the endpoints). 
Just for one thing, a polygon is not even a set of points on a plane but a sequence of such points;
see Proposition~\ref{prop:order} below and its proof.
Even the convexity of cyclic polygons is not a trivial matter; see \cite{jgeom,angles}.
Another cause of difficulties is that our definition of polygon convexity is broad enough not to require that any three vertices of a convex polygon be non-collinear, so that one also has to consider possible ``walks", possibly back and forth, along the 1-dimensional faces of the convex hull of the polygon.

In contrast, the proof of the ``upward'' theorem (Theorem \ref{th:conv-upwards}) is significantly
simpler and shorter than that of Theorem \ref{th:conv-reduction}. At first glance, this may now seem surprising, especially because Theorem~\ref{th:conv-upwards} may appear less intuitive that Theorem~\ref{th:conv-reduction}.

The paper is organized as follows. 
In Section \ref{results}, the definitions are given and the main results are stated: Theorems \ref{th:conv-reduction} and \ref{th:conv-upwards}; Propositions~\ref{prop:dim=1}, \ref{prop:order}, \ref{prop:G}, \ref{prop:def}, \ref{prop:apart}, \ref{prop:conv1}, and \ref{prop:conv-reduction}; and Corollaries~\ref{cor:sub-polygon}, \ref{cor:cut}, \ref{cor:test}, and \ref{cor:4-gons}. 
 
In Section \ref{proofs}, the necessary proofs are given. More specifically, Subsection \ref{lemmas, proofs of ths-props} contains statements of lemmas and based on them proofs of the main results stated in Section \ref{results}; the proofs of all lemmas are deferred further to Subsection \ref{proofs of lemmas}.  


\section{Definitions and results}\label{results}

A {\em polygon} is any finite sequence $\P:=(V_0,\dots,V_{n-1})$ of points (or, interchangeably, vectors) on the Euclidean plane $\R^2$. A polygon $\P=(V_0,\dots,V_{n-1})$, which is a sequence of $n$ points is also called an $n$-gon. The smallest value that one may allow for the integer $n$ is $0$, corresponding to a polygon with no vertices, that is, to the sequence $()$ of length $0$. 
The points $V_0,\dots,V_{n-1}$ are called the {\em vertices} of $\P$. 
The segments, or closed intervals,
$$[V_i,V_{i+1}]:=\conv\{V_i,V_{i+1}\}\quad\text{for}\ i\in\intr0{n-1}$$
are called the {\em edges} of polygon $\P$, where 
$$V_n:=V_0.$$
The symbol $\conv$ denotes, as usual, the convex hull \cite[page 12]{rock}. 
In general, our terminology corresponds to that in \cite{rock}. 
Here and in the sequel, we also use the notation
$$\intr km:=\{i\in\Z\colon k\le i\le m\},$$
where $\Z$ is the set of all integers; in particular, $\intr km$ is empty if $m<k$. 
Note that, if $V_i=V_{i+1}$, then the edge $[V_i,V_{i+1}]$ is a singleton set. 


Let us define the convex hull and dimension of polygon $\P$ as, respectively, the convex hull and dimension of the set of its vertices: $\conv\P:=\conv\{V_0,\dots,V_{n-1}\}$ and $\dim\P:=\dim\{V_0,\dots,V_{n-1}\}=\dim\conv\P$.


Now, a {\em convex polygon} can be defined as a polygon $\P$ such that the union of the edges of $\P$ coincides with the boundary $\partial\conv\P$ of the convex hull $\conv\P$ of $\P$; cf. e.g. \cite[page 5]{yaglom}. 
Thus, one has

\begin{definition}\label{def:conv}
A polygon $\P=(V_0,\dots,V_{n-1})$ is {\em convex} if 
$$\edg\P:=\bigcup_{i\in\intr0{n-1}}[V_i,V_{i+1}]=\partial\conv\P.$$
\end{definition}

One may note the following.

\begin{proposition}\label{prop:dim=1}
Any polygon $\P$ with $\dim\P\le1$ is convex.
\end{proposition}

\begin{remark}\label{rem:conv}
Let us emphasize that a polygon in this paper is a sequence and therefore ordered. In particular, even if the set $\{V_0,\dots,V_{n-1}\}$ of all vertices of a polygon $\P=(V_0,\dots,V_{n-1})$ coincides with the set $\ext\conv\P$ of all extreme points of the convex hull of $\P$, it does not necessarily follow that $\P$ is convex. For example, if $V_0=(0,0)$, $V_1=(1,0)$, $V_2=(1,1)$, and $V_3=(0,1)$, then polygon $(V_0,V_1,V_2,V_3)$ is convex, while polygon $(V_0,V_2,V_1,V_3)$ is not.
\end{remark}

\begin{definition}\label{def:strict}
A polygon $\P:=(V_0,\dots,V_{n-1})$ will be referred to as {\em strict} if the vertices $V_i$, $V_j$, and $V_k$ are non-collinear whenever $0\le i<j<k\le n-1$.
A polygon will be referred to as {\em strictly convex} if it is both strict and convex. 
\end{definition}

We shall now consider some other definitions of polygon convexity.

Let us recall \cite[page 31]{grun} that a {\em $2$-polytope} $\Pi$ is a compact closed set of dimension~2 in $\R^2$ such that the set $\ext\Pi$ of all extreme points of $\Pi$ is finite. 

\begin{proposition}\label{prop:order}
For any set $F\subseteq\R^2$, the following statements are equivalent to each other:
\begin{description}
\item[(i)]
$F=\ext\Pi$ for some $2$-polytope $\Pi$;
\item[(ii)]
$F=\{V_0,\dots,V_{n-1}\}$ for some strictly convex polygon $\P=(V_0,\dots,V_{n-1})$ with $n\ge3$.
\end{description}
\end{proposition}

As stated in Proposition \ref{prop:G} below, we shall show that a vertex enumeration 
$$\intr0{n-1}\ni i\longleftrightarrow V_i\in F$$
provided by Proposition~\ref{prop:order} is unique up to a cyclic permutation and/or an ``orientation switch". 

Let $\Xi_n$ denote the group of all permutations (i.e., bijections) of the set $\intr0{n-1}$ endowed with the composition operation. Let us denote the action and composition of the permutations in accordance with the formulas 
$$i\si:=\si(i)\quad\text{and}\quad i(\si\tau):=i\si\tau:=(i\si)\tau$$
for all $\si$ and $\tau$ in $\Xi_n$ and all $i\in\intr0{n-1}$.  

Let $G_n$ denote the subgroup of $\Xi_n$ generated by the ``primitive" cyclic permutation $\th$ and reflection $\rho$ defined by
$$
\begin{aligned}
(0\th,\dots,(n-1)\th)&:=(1,\dots,n-1,0)\quad\text{and} \\
(0\rho,\dots,(n-1)\rho)&:=(0,n-1,n-2,\dots,1).
\end{aligned}
$$
Since $\th\rho=\rho\th^{-1}$ and $\th^n=\rho^2=\text{the identity permutation}$, one has 
$$G_n= \{\th^j\rho^i\colon j\in\intr0{n-1},i\in\{0,1\}\}= \{\rho^i\th^j\colon j\in\intr0{n-1},i\in\{0,1\}\}.$$

For an $n$-gon $\P=(V_0,\dots,V_{n-1})$ and $\si\in\Xi_n$, define the corresponding permutation of $\P$ by the formula
$$\P\si:=(V_{0\si},\dots,V_{(n-1)\si}).$$

\begin{proposition}\label{prop:G}
Suppose that an $n$-gon $\P$ is strictly convex and $\si\in\Xi_n$. Then the $n$-gon $\P\si$ is strictly convex if and only if $\si\in G_n$.  
\end{proposition}

\begin{definition}\label{def:weak conv}
A polygon $\P$ will be referred to as {\em quasi-convex} if $$\edg\P\subseteq\partial\conv\P.$$
\end{definition}

Obviously, any convex polygon is quasi-convex.

\begin{definition}\label{def:to one side}
Let $P_0,\dots,P_m$ be any points on the plane, any two of which may in general coincide with each other. Let us write $P_2,\dots,P_m\|[P_0,P_1]$ 
and say
that
points $P_2,\dots,P_m$ are {\em to one side} 
of segment $[P_0,P_1]$ if 
there is a (straight) line $\ell$ containing $[P_0,P_1]$ and supporting to the set $\{P_0,\dots,P_m\}$; the latter, ``supporting" condition means here (in accordance with \cite[page 100]{rock}) that $\ell$ is the boundary of a closed 
half-plane containing the set $\{P_0,\dots,P_m\}$. 
Let us write $P_2,\dots,P_m\notOneSide[P_0,P_1]$ if it is not true that
$P_2,\dots,P_m\|[P_0,P_1]$. 
Let us say that
a polygon $\P=(V_0,\dots,V_{n-1})$ is to one side 
of a segment $[P_0,P_1]$ if the points $V_0,\dots,V_{n-1}$
are so.
Let us say that a polygon is {\em to-one-side} if it is to one side of every one of its edges. 
\end{definition}

\begin{proposition}\label{prop:def}
A polygon is quasi-convex if and only if 
it is to-one-side. 
\end{proposition}

The following proposition complements Proposition \ref{prop:def}.

\begin{proposition}\label{prop:apart}
Suppose that an $n$-gon $\P=(V_0,\dots,V_{n-1})$ is strictly convex and $\al$, $i$, and $\beta$ are integers such that
$1\le\al<i<\beta\le n-1$. Then $V_\al,V_\beta\notOneSide[V_0,V_i]$. 
\end{proposition}

\begin{definition}\label{def:ordinary}
A polygon $\P=(V_0,\dots,V_{n-1})$ is {\em ordinary} if 
its vertices are all distinct from one another: ($i\ne j\ \&\ i\in\intr0{n-1}\ \&\ j\in\intr0{n-1}\,)\implies V_i\ne V_j$.
\end{definition}

Note that in Proposition \ref{prop:def} the polygon is not assumed to be ordinary.

\begin{remark}\label{rem:ordinary}
The set of edges of any polygon can be represented as the union of the sets of edges of ordinary polygons.
As follows from \cite{grigoriev}, it takes $\Omega(n\log n)$ operations to test whether a polygon is ordinary.
\end{remark}


\begin{proposition}\label{prop:conv1}
An ordinary polygon is convex if and only if it is quasi-convex.
\end{proposition}

\begin{remark}\label{rem:conv1}
In Proposition \ref{prop:conv1}, one cannot drop the condition that the polygon is ordinary. For example, consider the polygon $(V_0,V_1,V_2,V_1)$ (with the edges $[V_0,V_1]$, $[V_1,V_2]$, $[V_2,V_1]$, $[V_1,V_0]$), where $V_0$, $V_1$, and $V_2$ are any three non-collinear points. Then $\P$ is quasi-convex but not convex, because $\frac12(V_0+V_2)\subseteq(\partial\conv\P)\setminus\edg\P$. 
\end{remark}

The following theorem is one of the main results of this paper.

\begin{theorem}\label{th:conv-reduction}
An ordinary convex polygon remains so after the removal of any one of its vertices. That is, if $\P=(V_0,\dots,V_{n-1})$ is an ordinary convex polygon, then the sub-polygon $\P^{(\al)}:=(V_0,\dots,V_{\al-1},V_{\al+1},\dots,V_{n-1})$ is so, for each $\al\in\intr0{n-1}$. 
\end{theorem}

\begin{remark}\label{rem:conv-reduction}
The condition in Theorem \ref{th:conv-reduction} that the polygon be ordinary cannot be dropped. For example, consider polygon $\P=(V_0,V_1,V_2,V_3,V_0,V_1,V_2,V_3)$, where $V_0=(0,0)$, $V_1=(1,0)$, $V_2=(1,1)$, and $V_3=(0,1)$, and let $\al=1$. Then polygon $\P$ is convex, while the sub-polygon $\P^{(\al)}=(V_0,V_2,V_3,V_0,V_1,V_2,V_3)$ is not convex. 
\end{remark}

If $\P=(V_0,\dots,V_{n-1})$ is a polygon, let us refer to any subsequence $(V_{i_0},\dots,V_{i_{m-1}})$ of $\P$, with $0\le i_0<\dots<i_{m-1}\le n-1$,  as a {\em sub-polygon} or, more specifically, 
as a {\em sub-$m$-gon} of $\P$.

Note that, according to this definition, a polygon of the form $\mathcal{Q}:=(V_{i_0},\dots,V_{i_{m-1}})$ 
with $0\le i_k<\dots<i_{m-1}<i_0<\dots<i_{k-1}\le n-1$
is generally not a sub-polygon of $\P=(V_0,\dots,V_{n-1})$; however, $\mathcal{Q}$ is a sub-polygon of the cyclic permutation $(V_{i_0},\dots,V_{n-1},V_0,\dots,V_{i_0-1})$ of polygon $\P$.

\begin{corollary}\label{cor:sub-polygon}
An ordinary convex polygon remains so after the removal of any number of its vertices. That is, if $\P=(V_0,\dots,V_{n-1})$ is an ordinary convex polygon, then any sub-polygon of $\P$ is so. 
\end{corollary}

Corollary~\ref{cor:sub-polygon} follows immediately from Theorem \ref{th:conv-reduction}.

\begin{corollary}\label{cor:cut}
Let $\P=(V_0,\dots,V_{n-1})$ be a convex polygon with $\Pi:=\conv\P$, and let $\ell$ be any line on the plane such that 
$\ell\cap\interi\Pi\ne\emptyset$, where $\interi$ denotes the interior. 
Then 
\begin{description}
	\item[(i)] 
$\ell\cap\edg\P=\{P,Q\}$, where $P$ and $Q$ are two distinct points such that $P\in[V_i,V_{i+1})$ and $Q\in[V_j,V_{j+1})$ for some $i$ and $j$ in $\intr0{n-1}$ with $i<j$ (here we use the notation $[A,B):=\{(1-t)A+tB\colon 0\le t<1\}$ if $A\ne B$, and $[A,B):=\emptyset$ if $A=B$);
	\item[(ii)] 
	if, moreover, the polygon $\P$ is strict then
the $i$ and $j$ in part (i) are uniquely determined, and the polygons $\P_+:=(V_0,\dots,V_i,P,Q,V_{j+1},\dots,V_{n-1})$ and $\P_-:=(P,V_{i+1},\dots,V_j,Q)$ are convex and lie to the opposite sides of the segment $[P,Q]$, in the sense that $\P$ is not to one side $[P,Q]$.	
\end{description}
\end{corollary}

Informally, Corollary~\ref{cor:cut} states that any straight line passing through the interior of the convex hull of a strictly convex polygon cuts it into two convex polygons; cf.\ Proposition~\ref{prop:apart}. 

\begin{remark}\label{rem:cut}
For the 
conclusion in part (ii) of Corollary~\ref{cor:cut} about the convexity of $\P_+$ and $\P_-$ to hold in general, the condition
that the polygon $\P$ be strict cannot be dropped; nor can it be relaxed to the condition that $\P$ be ordinary. Indeed, let $\P:=(V_0,V_1,V_2,V_3,V_4)$, where the points $V_0$, $V_1$, and $V_2$ are non-collinear, $V_3\in(V_0,V_2)$, and $V_4\in(V_2,V_3)$. 
Let also $\ell:=PQ$, where $P:=V_1$ and $Q\in(V_3,V_4)$.
Then it is not difficult to see that polygon $\P$ is ordinary.  


\vspace*{-35pt}

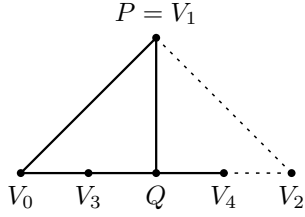
\parbox{4.1cm}
{
\pspicture(0,-2.2)(8,4)
\psset{unit=4.5mm}
\psset{origin={0,0}}
\psset{linewidth=.2mm}
\qdisk(0, 0){1.5pt}
\uput[270](0,0){$V_0$}
\qdisk(8, 0){1.5pt}
\uput[90](4,4){$P=V_1$}
\qdisk(4, 4){1.5pt}
\uput[270](8,0){$V_2$}
\qdisk(4, 0){1.5pt}
\uput[270](2,0){$V_3$}
\qdisk(2, 0){1.5pt}
\uput[270](6,0){$V_4$}
\qdisk(6, 0){1.5pt}
\uput[270](4,0){$Q$}
\psline[linewidth=.3mm](0,0)(6,0)
\psline[linewidth=.3mm](0,0)(4,4)
\psline[linewidth=.3mm](4,0)(4,4)
\psline[linestyle=dashed,dash=.1 .2](6,0)(8,0)
\psline[linestyle=dashed,dash=.1 .2](4,4)(8,0)
\endpspicture
}
\parbox{8cm}
{
On the other hand, in part (i) of Corollary~\ref{cor:cut} one must have $i=1$ and $j=3$, and polygon $\P_+=(V_0,V_1,V_1,Q,V_4)$ is not convex (nor is it quasi-convex), since $[V_1,Q]\not\subseteq\partial\conv\P_+$.
In this picture, the set $\edg\P_+$ consists of the solid line segments. 
}

\vspace*{-1.5cm}

\end{remark}

Theorem \ref{th:conv-reduction} is complemented by

\begin{proposition}\label{prop:conv-reduction}
Let $\P=(V_0,\dots,V_{n-1})$ be an ordinary convex polygon. 
If $V_\al\notin\ext\conv\P$ for some $\al\in\intr0{n-1}$, then the convex hull and union of the edges of the reduced polygon $\P^{(\al)}$ are the same as those of polygon $\P$. 
\end{proposition}

While Theorem \ref{th:conv-reduction} describes a downward hereditariness property of polygon convexity, the following theorem states that polygon convexity is also upward-hereditary. 

\begin{theorem}\label{th:conv-upwards}
Let $\P=(V_0,\dots,V_{n-1})$ be a polygon, with $n\ge5$. If the reduced polygon $\P^{(\al)}$ is convex for each $\al\in\intr0{n-1}$, then $\P$ is convex.
\end{theorem}

\begin{remark}\label{rem:n-upwards}
The condition $n\ge5$ in Theorem \ref{th:conv-upwards} cannot be dropped. Indeed, Theorem \ref{th:conv-upwards} cannot be true for $n=4$, because all theWhile Theorem \ref{th:conv-reduction} describes a downward hereditariness property of polygon convexity, the following theorem states that polygon convexity is also upward-hereditary. 
 $3$-gons are convex while not all $4$-gons are. On the other hand, Theorem \ref{th:conv-upwards} is trivially true for $n\in\{1,2,3\}$, because
all $n$-gons with $n\in\{1,2,3\}$ are convex.
\end{remark}

\begin{corollary}\label{cor:upwards}
Let $\P$ be an $n$-gon with $n\ge4$. Fix any $m\in\intr4n$. If all sub-$m$-gons of $\P$ are convex, then $\P$ is convex. 
\end{corollary}

This follows easily from Theorem~\ref{th:conv-upwards} by induction.

In particular, one has 

\begin{corollary}\label{cor:4-upwards} 
If all sub-$4$-gons of a polygon $\P$ are convex, then $\P$ is convex. 
\end{corollary}

\begin{corollary}\label{cor:test}
An ordinary $n$-gon $\P$ with $n\ge5$ is convex if and only if the reduced polygons $\P^{(\al)}$ are convex for all $\al\in\intr0{n-1}$. 
\end{corollary}

This follows immediately from Theorems \ref{th:conv-reduction} and \ref{th:conv-upwards}.

\begin{corollary}\label{cor:m-gons}
Let $\P$ be an ordinary $n$-gon with $n\ge4$.
Fix any $m\in\intr4n$. 
Then $\P$ is convex if and only if all sub-$m$-gons of $\P$ are convex. 
\end{corollary}

This follows immediately from Corollaries~\ref{cor:sub-polygon} and \ref{cor:upwards}. 

In particular, one has 

\begin{corollary}\label{cor:4-gons} 
An ordinary polygon $\P$ is convex if and only if all sub-$4$-gons of $\P$ are convex. 
\end{corollary}

Note that in Theorem~\ref{th:conv-upwards} and Corollaries~\ref{cor:upwards} and \ref{cor:4-upwards}
(in contrast with Theorem~\ref{th:conv-reduction} and Corollaries~\ref{cor:sub-polygon}, \ref{cor:test}, \ref{cor:m-gons}, and \ref{cor:4-gons}) it is not required that polygon $\P$ be ordinary.

A statement somewhat similar to 
Corollary~\ref{cor:4-upwards} was made in \cite{erdos} and reproduced in \cite{hall} as follows:
\begin{quote}
\textsc{Lemma 6.2.2.} {\em
If all the quadrilaterals formed from $n$ points, no three on a line, are convex, then the $n$ points are the vertices of a convex $n$-gon.\/}
\end{quote}
However, no explicit definition of the notion of a convex polygon (or that of a polygon in general) is found in \cite{erdos} or \cite{hall};
no proof of the quoted statement is given in \cite{erdos}, and the proof given in \cite{hall} is rather heuristic. From the context, it is apparent that the conclusion ``the $n$ points are the vertices of a convex $n$-gon" in the above quote should be interpreted as ``the set of the $n$ points is the set of all extreme points of a 2-polytope"; then it follows from Remark \ref{rem:conv} that this conclusion does not necessarily imply that the polygon (construed as a sequence) is convex according to Definition \ref{def:conv}; recall also Proposition~\ref{prop:order}. Note also that in Corollary~\ref{cor:4-gons} (in contrast with Proposition~\ref{prop:order})
it is not assumed that the polygon is strict.

For any two points $P$ and $Q$, let $(P,Q)$ denote the relative interior, $\ri[P,Q]$, of the segment $[P,Q]$, so that 
$(P,Q)=\ri[P,Q]=\{(1-t)P+tQ\colon 0<t<1\}$ if $P\ne Q$ and $(P,Q)=\ri[P,Q]=\emptyset$ otherwise. 
It is hoped that, within any given context, this notation will not be confused with that for the pair of points.

For any two distinct $P$ and $Q$, let $PQ$ denote the unique line containing the two points.

\section{Proofs}\label{proofs}

\subsection{Statements of lemmas, and proofs of the theorems and propositions}\label{lemmas, proofs of ths-props}

\begin{lemma}\label{lem:R.5}
Suppose that $k$ is a natural number and points $V_0,\dots,V_k$ are such that $\forall q\in\intr1{k-1}\ V_q\in[V_0,V_k]$. 
Then 
$$\bigcup_{q\in\intr0{k-1}}[V_q,V_{q+1}]=[V_0,V_k].$$
\end{lemma}

\begin{lemma}\label{lem:dim=1}
Proposition~\ref{prop:dim=1} is true: any polygon $\P$ with $\dim\P\le1$ is convex, and hence quasi-convex.
\end{lemma}

\begin{lemma}\label{lem:prop:def}
Proposition \ref{prop:def} is true: 
a polygon is quasi-convex if and only if 
it is to-one-side.
\end{lemma}

One reason for us to repeat here the statements of Propositions~\ref{prop:dim=1} and \ref{prop:def} is that, on the one hand, the facts stated in these
propositions will be used in the proofs of some other lemmas and, on the other hand, we want all the statements of the lemmas to precede the proofs of all propositions.

\begin{lemma}\label{lem:R.2}
Three points $P_0$, $P_1$, and $P_2$ on the plane are non-collinear if and only if for any vector $\nn$ on the plane one has the implication
$$\nn\cdot\p01=\nn\cdot\p02=0\implies\nn=\zero.$$
\end{lemma}

Here and elsewhere, as usual, the dot product is defined as the sum of the products of the respective coordinates.
\begin{lemma}\label{lem:R.1}
Let $P_0$, $P_1$, and $P_2$ be non-collinear points and let $\aa$, $\n{1}$, and $\n{2}$ be vectors (on the plane) such that $\n{1}\ne\zero$, $\n{2}\ne\zero$, and
$$\n{1}\cdot\p01=\n{2}\cdot\p02=0=\n{1}\cdot\aa=\n{2}\cdot\aa.$$
Then $\aa=\zero$.  
\end{lemma}

\begin{lemma}\label{lem:R.2a}
Suppose that a polygon $\P=(V_0,\dots,V_{n-1})$ is quasi-convex and $\ell$ is a supporting line to $\conv\P$. 
Then
$$\ell\cap\conv\P=[V_i,V_k]$$ 
for some $i$ and $k$ in $\intr0{n-1}$, and, moreover, $V_i$ and $V_k$ are extreme points of $\conv\P$. 
\end{lemma}

\begin{lemma}\label{lem:R.3}
Suppose that one has the following conditions: (i) a polygon $\P=(V_0,\dots,V_{n-1})$ is quasi-convex; (ii) $\ell$ is a supporting line to $\conv\P$; (iii) $j\in\intr1n$; (iv) $V_{j-1}\in\ell$; and (v) $V_j\notin\ell$. Then $V_{j-1}\in\ext\conv\P$.
\end{lemma}

\begin{lemma}\label{lem:extreme}
If $\P=(V_0,\dots,V_{n-1})$ is an ordinary quasi-convex polygon of dimension $\dim\P=2$ and $V_1\in\ext\conv\P$, then the points $V_0$, $V_1$, and $V_2$ are non-collinear. 
\end{lemma}

\begin{lemma}\label{lem:case Val extreme}
Suppose that $\P=(V_0,\dots,V_{n-1})$ is an ordinary quasi-convex polygon and $V_\al\in\ext\conv\P$, for some $\al\in\intr0{n-1}$.
Then $\P^{(\al)}=(V_0,\dots,V_{\al-1},V_{\al+1},\dots,\\
V_{n-1})$ is also an ordinary quasi-convex polygon. 
\end{lemma}

\begin{lemma}\label{lem:R.13}
Suppose that a polygon $\P=(V_0,\dots,V_{n-1})$ is ordinary and quasi-convex, a line $\ell$ is supporting to $\conv\P$, and $\ell\cap\conv\P=[V_0,V_k]$ for some $k\in\intr0{n-1}$. Then one has the implication
$$\bigl(\exists q\in\intr1{k-1}\ V_q\notin\ell\bigr)
\implies\bigl(\forall q\in\intr1{k-1}\ V_q\notin\ell\bigr);$$
equivalently, one has the implication
$$\bigl(\exists q\in\intr1{k-1}\ V_q\in\ell\bigr)
\implies\bigl(\forall q\in\intr1{k-1}\ V_q\in\ell\bigr).$$
\end{lemma}

\begin{lemma}\label{lem:prop:conv-reduction}
Let $\P=(V_0,\dots,V_{n-1})$ be an ordinary quasi-convex polygon.
If $V_\al\notin\ext\conv\P$ for some $\al\in\intr0{n-1}$, then the convex hull and union of the edges of the reduced polygon $\P^{(\al)}=(V_0,\dots,V_{\al-1},V_{\al+1},\dots,V_{n-1})$ are the same as those of polygon $\P$:
$$\conv\P^{(\al)}=\conv\P\quad\text{and}\quad
\edg\P^{(\al)}=\edg\P,$$
so that $\P$ is an ordinary quasi-convex polygon.  
\end{lemma}

\begin{lemma}\label{lem:conv-reduction}
Suppose that $\P=(V_0,\dots,V_{n-1})$ is an ordinary quasi-convex polygon and $\al\in\intr0{n-1}$.
Then $\P^{(\al)}=(V_0,\dots,V_{\al-1},V_{\al+1},\dots,V_{n-1})$ is also an ordinary quasi-convex polygon. 
\end{lemma}

\begin{lemma}\label{lem:new}
Suppose that $n\ge5$ and $\P=(V_0,\dots,V_{n-1})$ is a quasi-convex polygon such that $\P^{(\al)}=(V_0,\dots,V_{\al-1},V_{\al+1},\dots,V_{n-1})$ is convex for every $\al\in\intr0{n-1}$.
Then polygon $\P$ is convex. 
\end{lemma}

\begin{lemma}\label{lem:R.10}
The following conditions can never take place all together:
$\P=(V_0,V_1,V_2,V_3)$ is an ordinary quasi-convex polygon; 
a line $\ell$ is supporting to $\conv\P$; 
$\{V_0,V_2\}\subseteq\ell$; and
$\{V_1,V_3\}\cap\ell=\emptyset$. 
\end{lemma}

\begin{lemma}\label{lem:R.11}
The following conditions can never take place all together:
$\P=(V_0,\dots,V_{n-1})$ is an ordinary quasi-convex polygon; 
a line $\ell$ is supporting to $\conv\P$; 
$\ell\cap\conv\P=[V_0,V_k]$ for some $k\in\intr2{n-2}$;
$\{V_1,\dots,V_{k-1}\}\cap\ell=\emptyset$; and \break
$\{V_{k+1},\dots,V_{n-1}\}\cap\ell=\emptyset$.  
\end{lemma}

\begin{lemma}\label{lem:conv1}
Proposition~\ref{prop:conv1} is true: an ordinary polygon is convex if and only if it is quasi-convex.   
\end{lemma}

\begin{lemma}\label{lem:th:conv-reduction}
Theorem~\ref{th:conv-reduction} is true: if $\P=(V_0,\dots,V_{n-1})$ is an ordinary convex polygon, then the sub-polygon $\P^{(\al)}=(V_0,\dots,V_{\al-1},V_{\al+1},\dots,V_{n-1})$ is so, for each $\al\in\intr0{n-1}$.   
\end{lemma}

\begin{lemma}\label{lem:cor:sub-polygon}
Corollary~\ref{cor:sub-polygon} is true: if $\P=(V_0,\dots,V_{n-1})$ is an ordinary convex polygon, then any sub-polygon of $\P$ is so.
\end{lemma}

\begin{lemma}\label{lem:sep}
Let $\Pi$ be a $2$-polytope, and let $V$ be any point in the set $\R^2\setminus\Pi$. 
Then there exists a line $\ell$ such that $\dim(\ell\cap\Pi)=1$ and $\ell=\partial H$ for some closed half-plane $H$ such that $V\notin H\supseteq\Pi$. 
\end{lemma}

\begin{lemma}\label{lem:Q}
Let $\Q:=(U_0,U_1,U_2,U_3)$ be such a $4$-gon that there exists a point $O\in[U_0,U_2]\cap(U_1,U_3)$. Let $\ell$ be a line containing point $U_1$ and supporting to the ``triangle" polytope $\conv\{U_0,U_1,U_2\}$. Then $\ell$ is supporting to polytope $\conv\Q$ as well.
\end{lemma}

\begin{lemma}\label{lem:prop:apart}
Proposition~\ref{prop:apart} is true: if an $n$-gon $\P=(V_0,\dots,V_{n-1})$ is strictly convex and $\al$, $i$, and $\beta$ are integers such that
$1\le\al<i<\beta\le n-1$, then $V_\al,V_\beta\notOneSide[V_0,V_i]$. 
\end{lemma}

\begin{lemma}\label{lem:if}
If an $n$-gon $\P$ is strictly convex and $\si\in G_n$, then the $n$-gon $\P\si$ is strictly convex as well.  
\end{lemma}

\begin{lemma}\label{lem:0213}
If a $4$-gon $\Q:=(U_0,U_1,U_2,U_3)$ is strictly convex, then the $4$-gon $\tilde\Q:=(U_0,U_2,U_1,U_3)$ is not.
\end{lemma}

\begin{lemma}\label{lem:si(1)}
Suppose that $\si\in\Xi_n$, $0\si=0$, $1\si<(n-1)\si$, and
$n$-gons $\P=(V_0,\dots,V_{n-1})$ and $\P\si$ are both strictly convex.
Then $1\si=1$. 
\end{lemma}

\begin{lemma}\label{lem:si}
Suppose that $\si\in\Xi_n$, $0\si=0$, $1\si<(n-1)\si$, and
$n$-gons $\P$ and $\P\si$ are both strictly convex.
Then $i\si=i$ for all $i\in\intr0{n-1}$; i.e., $\si$ is the identity permutation. 
\end{lemma}

\begin{lemma}\label{lem:al-beta}
If $\P=(V_0,\dots,V_{n-1})$ is a strictly convex polygon and integers $\al$ and $\beta$ in $\intr0{n-1}$ are distinct from each other, then 
$[V_\al,V_{\al+1})\cap[V_\beta,V_{\beta+1})=\emptyset$.
\end{lemma}

\begin{proof}[Proof of Proposition~\ref{prop:dim=1}]
This follows immediately from Lemma~\ref{lem:dim=1}.
\end{proof}

\begin{proof}[Proof of Proposition \ref{prop:order}]
{\bf (ii)$\implies$(i)}\quad Suppose that $F=\{V_0,\dots,V_{n-1}\}$ for some strictly convex polygon $\P=(V_0,\dots,V_{n-1})$ with $n\ge3$. Let 
$\Pi:=\conv F=\conv\P$. Then, by \cite[Corollary 18.3.1]{rock}, one has $\ext\Pi\subseteq F$. 
Let $\ell:=V_0V_1$. By Lemma~\ref{lem:prop:def}, the line $\ell$ is supporting to $\conv\P$. Since polygon $\P=(V_0,\dots,V_{n-1})$ is strict and $\ell=V_0V_1$, one has $V_2\notin\ell$.
Hence, by Lemma \ref{lem:R.3}, $V_1\in\ext\conv\P=\ext\Pi$. 
Similarly, $V_j\in\ext\Pi$ for all $j\in\intr0{n-1}$, so that $F\subseteq\ext\Pi$ and hence $F=\ext\Pi$.
That $\dim\Pi=2$ follows because $\Pi=\conv\P$ and $\P$ is a strict $n$-gon with $n\ge3$. 

{\bf (i)$\implies$(ii)}\quad Suppose that $F=\ext\Pi$ for some $2$-polytope $\Pi$, so that $n:=\card F$ is finite. If the implication (i)$\implies$(ii) fails for some natural $n$, then let $n$ be the smallest such number. Note that one must have $n\ge3$, since $\Pi$ is a $2$-polytope and hence $\dim\Pi=2$. 

Consider first the case $n=3$, so that $F=\{V_0,V_1,V_2\}$ for some non-collinear points $V_0,V_1,V_2$ in $\R^2$. Then polygon $\P:=(V_0,V_1,V_2)$ is trivially to-one-side. Hence, in view of Lemmas \ref{lem:prop:def} and \ref{lem:conv1}, polygon $\P$ is strictly convex.

Let now $n\ge4$. Since $F=\ext\Pi$ for some $2$-polytope $\Pi$, 
one has $\Pi=\conv F$, by \cite[Theorem 18.5]{rock}. Hence,
there exist three non-collinear points $U_0,U_1,U_2$ in $F$. Moreover, since $n\ge4$, there exists a point
$V\in F\setminus\{U_0,U_1,U_2\}$. Let now $\tilde F:=F\setminus\{V\}$ and $\tilde\Pi:=\conv\tilde F$.
Note that $\dim\tilde\Pi=\dim\tilde F=2$, because the set $\tilde F$ contains the three non-collinear points $U_0,U_1,U_2$.

Note also that $\ext\tilde\Pi=\tilde F$. 
$\bigl($Indeed, $\ext\tilde\Pi=\ext\conv\tilde F\subseteq\tilde F$, again by \cite[Corollary 18.3.1]{rock}. On the other hand, 
$\tilde F\subseteq\tilde\Pi\cap F=\tilde\Pi\cap\ext\Pi\subseteq\ext\tilde\Pi$.$\bigr)$ 

Hence and because of the minimality of $n$, there exists a strictly convex polygon $\tilde\P=(V_0,\dots,V_{n-2})$ such that $\tilde F=\{V_0,\dots,V_{n-2}\}$, whence $\tilde\Pi=\conv\tilde\P$. 
Observe that $V\notin\tilde\Pi$. $\Bigl($Otherwise, 
$V\in\tilde\Pi=\conv\{V_0,\dots,V_{n-2}\}$.
Hence, by \cite[Corollary 2.3.1]{rock}, there exist nonnegative real numbers $\la_0,\dots,\la_{n-2}$ such that $\la_0+\dots+\la_{n-2}=1$ and $\la_0V_0+\dots+\la_{n-2}V_{n-2}=V$. Since 
$V\notin\tilde F=\{V_0,\dots,V_{n-2}\}$, at least one of the $\la_i$'s (say $\la_0$) belongs to the interval $(0,1)$. Then, for $U:=(\la_1V_1+\dots+\la_{n-2}V_{n-2})/(1-\la_0)$, one has $V=\la_0V_0+(1-\la_0)U$, $V_0\in\Pi$, and $U\in\Pi$, which contradicts the conditions $V\in F=\ext\Pi$.$\Bigr)$

Therefore, in view of Lemma \ref{lem:sep}, there exists a line $\ell$ such that $\dim(\ell\cap\tilde\Pi)=1$ and $\ell=\partial H$ for some closed half-plane $H$ such that $V\notin H\supseteq\tilde\Pi$.
By Lemma \ref{lem:R.2a}, $\ell\cap\tilde\Pi=[V_i,V_k]$ for some $i$ and $k$ in $\intr0{n-2}$, and $V_i$ and $V_k$ are in $\ext\tilde\Pi$.
Hence,
$$[V_i,V_k]=\ell\cap\tilde\Pi
=\tilde\Pi\cap\partial H
\subseteq\partial\tilde\Pi
=[V_{n-2},V_0]\cup\bigcup_{q=0}^{n-3}[V_q,V_{q+1}]$$
(since $\tilde\Pi=\conv\tilde\P$ and $\tilde\P=(V_0,\dots,V_{n-2})$ is convex).
The conditions $\dim(\ell\cap\tilde\Pi)=1$ and $\ell\cap\tilde\Pi=[V_i,V_k]$ imply that $V_i\ne V_k$, so that the set $[V_i,V_k]$ is infinite and hence contains at least two distinct points in one of the intervals $[V_{n-2},V_0],[V_0,V_1],\dots,[V_{n-3},V_{n-2}]$, covering the set $[V_i,V_k]$.

W.l.o.g., these two distinct points lie in the interval $[V_{n-2},V_0]$, so that
$[V_{n-2},V_0]\subseteq\ell\cap\tilde\Pi=[V_i,V_k]$. 
Since $\{V_{n-2},V_0\}\subseteq\tilde F=\ext\tilde\Pi$, it follows that $\{V_{n-2},V_0\}\subseteq\ext(\ell\cap\tilde\Pi)$. 
Because polygon $\P$ is strict and $n\ge4$, one has $V_{n-2}\ne V_0$ and hence
$$\ell\cap\tilde\Pi=[V_{n-2},V_0].$$

Because $\Pi\supseteq\tilde\Pi$, it follows that $\ell\cap\Pi\supseteq[V_{n-2},V_0]$. 
On the other hand, if there existed a point $U\in\ell\cap\Pi\setminus[V_{n-2},V_0]$, then one would have either $V_0\in(U,V_{n-2})$ or $V_{n-2}\in(U,V_0)$, which would contradict the conditions $\{V_0,V_{n-2}\}\subseteq F=\ext\Pi$. 
One concludes that
$$\ell\cap\Pi=[V_{n-2},V_0].$$

Let now
$$V_{n-1}:=V$$
and
$$\P:=(V_0,\dots,V_{n-1}).$$
Since
$$\{V_0,\dots,V_{n-1}\}=F=\ext\Pi$$
and no three distinct extremal points of a convex set can be collinear, it follows that the polygon $\P$ is strict and hence ordinary. 

In view of Lemma~\ref{lem:conv1}, it remains to show that $\P$ is quasi-convex. To this end, recall first that $\ell=\partial H$, where $H$ is a close half-plane such that 
$V_{n-1}=V\notin H\supseteq\tilde\Pi$. 

Let $H_+$ be the interior of $H$ and $H_-:=\R^2\setminus H$, so that $H_+$ and $H_-$ are the two open half-planes whose common boundary is the line $\ell$. 
Then $V_{n-1}\in\R^2\setminus H=H_-$. Also, $V_i\in H\setminus\ell=H_+$ for all $i\in\intr1{n-3}$ (because $\{V_1,\dots,V_{n-3}\}\subseteq\tilde\Pi\subseteq H$ while $$\{V_1,\dots,V_{n-3}\}\cap\ell=\{V_1,\dots,V_{n-3}\}\cap\Pi\cap\ell
=\{V_1,\dots,V_{n-3}\}\cap[V_{n-2},V_0]=\emptyset,$$ 
since polygon $\P=(V_0,\dots,V_{n-1})$ is strict). 

The strictness of polygon $\P=(V_0,\dots,V_{n-1})$ also implies that, for each 
$$i\in\intr1{n-3},$$
the line 
$$\tilde\ell:=V_iV_{n-1}$$
differs from the line $\ell=V_0V_{n-2}$. 
Moreover, if line $\tilde\ell$ were parallel to $\ell$, then one would have $\tilde\ell=V_iV_{n-1}\subseteq H_-=\R^2\setminus H$ (since $\partial H_-=\ell$ and $V_{n-1}\in H_-$), which would in turn imply that $V_i\in\tilde\ell\subseteq\R^2\setminus H$, which would contradict the condition $V_i\in H$. 
Hence,
$$\ell\cap\tilde\ell=\{O\},$$
for some point $O$.  

Since $O\in\ell\subseteq\R^2\setminus H_+$ and 
$V_{n-1}\in H_-\subseteq\R^2\setminus H_+$, one has $[O,V_{n-1}]\subseteq\R^2\setminus H_+$. This and the condition $V_i\in H_+$ imply that $V_i\notin[O,V_{n-1}]$. 
Similarly, because $O\in\ell\subseteq\R^2\setminus H_-$ and $V_i\in H_+\subseteq\R^2\setminus H_-$, one has $[O,V_i]\subseteq\R^2\setminus H_-$. This and the condition $V_{n-1}\in H_-$ imply that $V_{n-1}\notin[O,V_i]$. 
Since $O\in\tilde\ell=V_iV_{n-1}$, it follows that 
$$O\in(V_i,V_{n-1}).$$

Note also that $O\in\ell\cap(V_i,V_{n-1})\subseteq\ell\cap\Pi=[V_0,V_{n-2}]$. 

In view of Lemma \ref{lem:prop:def}, the convexity of polygon $\tilde\P=(V_0,\dots,V_{n-2})$, and the strictness of polygon $\P=(V_0,\dots,V_{n-1})$, it suffices to show that
\begin{description}
\item[(i)] $V_{n-1},V_0,V_{n-2}\|[V_j,V_{j+1}]$ for all $j\in\intr0{n-3}$\quad and 
\item[(ii)] $V_i,V_0\|[V_{n-2},V_{n-1}]$ and $V_i,V_{n-2}\|[V_{n-1},V_0]$ for all $i\in\intr1{n-3}$.
\end{description}

Let now $\ell_j:=V_jV_{j+1}$, where $j\in\intr0{n-1}$. 

If $j\in\intr0{n-3}$, then there exists some $i\in\intr1{n-3}$ such that $V_i\in\ell_j$. Moreover, then the line $\ell_j$ is supporting to $\tilde\Pi=\conv\{V_0,\dots,V_{n-2}\}$ (by Lemma \ref{lem:prop:def}, since polygon $\tilde\P=(V_0,\dots,V_{n-2})$ is convex); hence, $\ell_j$ is supporting also to the ``triangle" polytope $\conv\{V_0,V_i,V_{n-2}\}$. 

Otherwise, if $j\in\{n-2,n-1\}$, then $V_{n-1}\in\ell_j$. Moreover, then the line $\ell_j$ (which is either $V_{n-2}V_{n-1}$ or $V_{n-1}V_0$) is obviously supporting to the ``triangle" polytope $\conv\{V_0,V_{n-2},V_{n-1}\}$. 

Therefore, it remains to apply Lemma \ref{lem:Q} with $\ell_j$ in place of $\ell$ (with $\Q=(U_0,U_1,U_2,U_3):=(V_0,V_i,V_{n-2},V_{n-1})$ to prove item~(i) above and with $\Q=(U_0,U_1,U_2,U_3):=(V_0,V_{n-1},V_{n-2},V_i)$ to prove item~(ii) ). 
Proposition~\ref{prop:order} is proved.
\end{proof}

\begin{proof}[Proof of Proposition \ref{prop:G}]
The ``if" part follows immediately from Lemma \ref{lem:if}. 

To prove the ``only if" part, suppose that indeed $n$-gons $\P$ and $\P\si$ are strictly convex for some $\si\in\Xi_n$. We have to show that then $\si\in G_n$. W.l.o.g., $n\ge3$, because $\Xi_n=G_n$ for $n\le2$. 

Letting $j:=0\si$, one has $0\la=0$ for $\la:=\si\th^{-j}$, where $\th$ is the ``primitive" cyclic permutation defined before the statement of Proposition \ref{prop:G}. Since $n\ge3$, one also has $1\ne n-1$ and hence $1\la\ne(n-1)\la$. 
Let $\mu:=\la$ if $1\la<(n-1)\la$ and $\mu:=\la\rho$ if, otherwise, $1\la>(n-1)\la$, where $\rho$ is the reflection permutation defined before the statement of Proposition \ref{prop:G}. 
Then $0\mu=0$ and $1\mu<(n-1)\mu$. 
On the other hand, in view of Lemma \ref{lem:if} and strict convexity of $\P\si$, polygon $\P\mu$ is strictly convex, because permutations $\th^{-j}$ and $\rho$ belong to the group $G_n$. 
Now Lemma \ref{lem:si} implies that $\mu$ is the identity permutation. Thus, the permutation $\si$ equals either $\th^j$ or $\rho^{-1}\th^j$, so that $\si\in G_n$. 
\end{proof}

\begin{proof}[Proof of Proposition \ref{prop:def}]
This follows immediately from Lemma \ref{lem:prop:def}.
\end{proof}

\begin{proof}[Proof of Proposition \ref{prop:apart}]
This follows immediately from Lemma~\ref{lem:prop:apart}
\end{proof}

\begin{proof}[Proof of Proposition \ref{prop:conv1}]
This follows immediately from Lemma \ref{lem:conv1}.
\end{proof}

\begin{proof}[Proof of Theorem \ref{th:conv-reduction}]
This follows immediately from Lemma~\ref{lem:th:conv-reduction}.
\end{proof}

\begin{proof}[Proof of Corollary~\ref{cor:cut}]
\textbf{(i)}\quad The set $\ell\cap\Pi$ is compact and convex, and it has a nonempty relative interior (since $\ell\cap\interi\Pi\ne\emptyset$), so that $\ell\cap\Pi=[P,Q]$ for some distinct points $P$ and $Q$, which necessarily lie on the boundary $\partial\Pi$. 
On the other hand, by \cite[Theorem~6.1]{rock}, $(P,Q)\subseteq\interi\Pi$. 
Therefore and because polygon $\P$ is convex, $\ell\cap\edg\P=\ell\cap\partial\conv\P=\ell\cap\partial\Pi=\{P,Q\}$, so that $\{P,Q\}\subseteq\edg\P=\bigcup_{i=0}^{n-1}[V_i,V_{i+1})$.
Thus, $P\in[V_i,V_{i+1})$ and $Q\in[V_j,V_{j+1})$ for some $i$ and $j$ in $\intr0{n-1}$, and w.l.o.g.\ $i\le j$. If at that $i=j$, then the two distinct points, $P$ and $Q$, lie in $[V_i,V_{i+1})$, so that $[V_i,V_{i+1})$ is an infinite subset of the line $PQ=\ell$, which contradicts the above conclusion that $\ell\cap\edg\P=\{P,Q\}$. It follows that $i<j$, which completes the proof of part (i) of Corollary~\ref{cor:cut}.

\textbf{(ii)}\quad 
Here it is assumed that $\P$ is strict, and hence ordinary. Then, in view of Lemma~\ref{lem:al-beta}, the $i$ and $j$ in part (i) are uniquely determined.

To prove the rest of part (ii) of Corollary~\ref{cor:cut}, let us first consider the case when $P\ne V_i$ and $Q\ne V_j$.
Then define the ``extended'' polygon by the formula $$\hat\P:=(V_0,\dots,V_i,P,V_{i+1},\dots,V_j,Q,V_{j+1},\dots,V_{n-1}).$$ 

Note that polygon $\hat\P$ is ordinary.
Indeed, the vertices $V_0,\dots,V_{n-1}$ of $\P$ are distinct from one another. Also, it was seen that $P\ne Q$, and the assumptions $P\ne V_i$ and $Q\ne V_j$ imply that $P\in(V_i,V_{i+1})$ and $Q\in(V_j,V_{j+1})$, so that  $\{P,Q\}\cap\{V_0,\dots,V_{n-1}\}
=\emptyset$ (since polygon $\P$ is strict). 

Besides, $\conv\hat\P=\conv\P$ and $\edg\hat\P=\edg\P$, so that polygon $\hat\P$ inherits the convexity property of $\P$. 
The polygons $\P_\pm$ are sub-polygons of $\hat\P$ and hence convex, by Lemma~\ref{lem:cor:sub-polygon}.

Now notice that $\P$ is not to one side of $[P,Q]$. Indeed, otherwise the line $\ell=PQ$ would be supporting to $\Pi$, which would contradict the condition $\ell\cap\interi\Pi\ne\emptyset$.

This completes the proof of part (ii) in the case when $P\ne V_i$ and $Q\ne V_j$. The other three cases, depending on which of the equalities $P=V_i$ and $Q=V_j$ hold(s), are considered quite similarly, with the definition of the extended polynomial appropriately modified. 

For instance, in the case when $P=V_i$ and $Q\ne V_j$,
define the ``extended'' polygon by the formula $$\hat\P:=(V_0,\dots,V_i,V_{i+1},\dots,V_j,Q,V_{j+1},\dots,V_{n-1}),$$ 
so that it be ordinary. Then (in contrast with $\P_-$) the polygon $\P_+$ is generally \emph{not} a sub-polygon of $\hat\P$. However, the polygon $\tilde\P_+:=(V_0,\dots,V_i,Q,V_{j+1},\dots,V_{n-1})$ is so and hence is convex. 
On the other hand, $\conv\P_+=\conv\tilde\P_+$ and $\edg\P_+=\edg\tilde\P_+$, so that $\P_+$ is also convex.

This completes the proof of entire Corollary~\ref{cor:cut}.
\end{proof}

\begin{proof}[Proof of Proposition \ref{prop:conv-reduction}]
This follows immediately from 
Lemma \ref{lem:conv1} and Lemma \ref{lem:prop:conv-reduction}.
\end{proof}

\begin{proof}[Proof of Theorem \ref{th:conv-upwards}]
Suppose that $\P=(V_0,\dots,V_{n-1})$, $n\ge5$, and $\P^{(\al)}$ is convex for each $\al\in\intr0{n-1}$.
By Lemma~\ref{lem:new}, it suffices to show that 
$\P$ is quasi-convex.
Assume the contrary. 
Then, by Lemma \ref{lem:prop:def}, there exist $\beta$, $j$, and $k$ in $\intr0{n-1}$ such that $V_j,V_k\notOneSide[V_\beta,V_{\beta+1}]$.
Without loss of generality (w.l.o.g.), $\beta=0$ (otherwise, consider the cyclic permutation
$(W_0,\dots,W_{n-1}):=(V_\beta,\dots,V_{n-1},V_0,\dots,V_{\beta-1})$ of the vertices). 
Since $n\ge5$, there exists some $\al\in\intr0{n-1}\setminus\{0,1,j,k\}$, so that $\{0,1,j,k\}\subseteq\intr0{n-1}\setminus\{\al\}$ and $\al\in\intr2{n-1}$. 
For $q\in\intr0{n-2}$, let
$$U_q:=V_{\hat q},$$
where
$$\hat q:=
\begin{cases}
q &\text{if}\ q\in\intr0{\al-1},\\
q+1 &\text{if}\ q\in\intr\al{n-2}.
\end{cases}
$$
Then $\P^{(\al)}=(U_0,\dots,U_{n-2})$ and
$$Q:=\{\hat q\colon q\in\intr0{n-2}\}=\intr0{n-1}\setminus\{\al\}
\supseteq\{0,1,j,k\}.$$
Hence, there exists $q$ and $r$ in $\intr0{n-2}$ such that $\hat q=j$ and $\hat r=k$, whence $V_j=U_q$ and $V_k=U_r$. 
Also, the condition $\al\in\intr2{n-1}$ implies that $\hat q=q$ for $q\in\{0,1\}$, so that $U_0=V_0$ and $U_1=V_1$. 
Recall that the polygon $\P^{(\al)}=(U_0,\dots,U_{n-2})$ is convex. Hence, by Lemma \ref{lem:prop:def}, one has $U_q,U_r\|[U_0,U_1]$; that is, $V_j,V_k\|[V_0,V_1]$, which contradicts the assumption $V_j,V_k\notOneSide[V_\beta,V_{\beta+1}]$ (with $\beta=0$). 
\end{proof}

\subsection{Proofs of the lemmas}\label{proofs of lemmas}

\begin{proof}[Proof of Lemma \ref{lem:R.5}]
Let $\ell$ be a line. 
Observe that, if $\de_1\subseteq\ell$ and $\de_2\subseteq\ell$ are intervals such that $\de_1\cap\de_2\ne\emptyset$, then $\de_1\cup\de_2$ is also an interval.
Observe also that for all $j\in\intr1{k-1}$ the set 
$$\lp\bigcup_{q\in\intr0{j-1}}[V_q,V_{q+1}]\rp
\cap[V_j,V_{j+1}]$$
is non-empty, since it contains point $V_j$.
Now it follows
by induction in $k$ that $\bigcup_{q=0}^{k-1}[V_q,V_{q+1}]$
is an interval. Since this interval contains the points $V_0$ and $V_k$, one has $\bigcup_{q=0}^{k-1}[V_q,V_{q+1}]\supseteq[V_0,V_k]$. 
On the other hand, $\forall q\in\intr1{k-1}\ V_q\in[V_0,V_k]$ and hence $\bigcup_{q=0}^{k-1}[V_q,V_{q+1}]\subseteq[V_0,V_k]$.
\end{proof}

\begin{proof}[Proof of Lemma~\ref{lem:dim=1}]
This follows easily from Lemma~\ref{lem:R.5}.
Indeed, let $\P=(V_0,\dots,\\
V_{n-1})$ be a polygon with $\dim\P\le1$. Then there is a line $\ell$ containing all the vertices $V_0,\dots,V_{n-1}$. Consider the order on $\ell$ induced by any one-to-one affine mapping of $\R$ onto $\ell$. Let $i$ and $k$ in $\intr0{n-1}$ be such that $V_i$ and $V_k$ are, respectively, the minimum and the maximum of the set of $n$ vertices of polygon $\P$, according to the chosen order. W.l.o.g., $i=0$, so that $V_q\in[V_0,V_k]$ $\forall q\in\intr0{n-1}$, whence $\partial\conv\P=\conv\P=[V_0,V_k]$. Also, it follows by Lemma~\ref{lem:R.5} that
$\edg\P\supseteq\bigcup\limits_{q\in\intr0{k-1}}[V_q,V_{q+1}]=[V_0,V_k]\supseteq\edg\P$, whence
$\edg\P=[V_0,V_k]=\partial\conv\P$. 
\end{proof}

\begin{proof}[Proof of Lemma \ref{lem:prop:def}]

{\bf ``Only if":}\quad Suppose that a polygon $\P=(V_0,\dots,V_{n-1})$ is quasi-convex. The case $\dim\P\le1$ is easy, because then there is a line $\ell$ containing the entire convex hull $\conv\P$, so that $\ell$ contains all the edges of polygon $\P$ and is supporting to $\conv\P$. Let now $\dim\P=2$. Take any $i\in\intr0{n-1}$. Then, by Definition \ref{def:weak conv}, $[V_i,V_{i+1}]\subseteq\partial\conv\P$. Hence, by \cite[Theorem 11.6]{rock}, there is a line $\ell$ containing $[V_i,V_{i+1}]$ and supporting to $\conv\P$. 

{\bf ``If":}\quad Suppose a polygon $\P=(V_0,\dots,V_{n-1})$ is to-one-side. Take any $i\in\intr0{n-1}$. Then there is a line $\ell$ containing $[V_i,V_{i+1}]$ and supporting to $\conv\P$. Let $H$ be the corresponding half-plane containing $\conv\P$. By the definition of the convex hull, $[V_i,V_{i+1}]\subseteq\conv\P$. Hence,
$$[V_i,V_{i+1}]\subseteq\ell\cap\conv\P=\partial H\cap\conv\P\subseteq\partial \conv\P$$
(the latter inclusion follows because $\conv\P\subseteq H$). 
\end{proof}

\begin{proof}[Proof of Lemma \ref{lem:R.2}]
{\bf ``Only if"}\quad Assume that $\nn\cdot\p01=\nn\cdot\p02=0$ while $\nn\ne\zero$. Then the points $P_0$, $P_1$, and $P_2$ lie on the line $\ell:=\{P\colon\nn\cdot\overrightarrow{P_0P}=0\}$, which is a contradiction. 

{\bf ``If"}\quad Assume that points $P_0$, $P_1$, and $P_2$ are collinear, so that they lie on one line, which must be a set of the form $\{P\colon\nn\cdot P=c\}$ for some vector $\nn\ne\zero$ and some real number $c$. Then $\nn\cdot\p01=\nn\cdot\p02=0$ while $\nn\ne\zero$, which is a contradiction.
\end{proof}

\begin{proof}[Proof of Lemma \ref{lem:R.1}]
Assume the contrary: that $P_0$, $P_1$, and $P_2$ are non-collinear points, $\n{1}\ne\zero$, $\n{2}\ne\zero$, 
$\n{1}\cdot\p01=\n{2}\cdot\p02=0=\n{1}\cdot\aa=\n{2}\cdot\aa$, while $\aa\ne\zero$. 
Then points $\zero$, $\n{1}$, and $\n{2}$ are collinear, since they all lie on the line $\{\nn\in\R^2\colon\aa\cdot\nn=0\}$. Since $\n{1}\ne\zero$ and $\n{2}\ne\zero$, one has $\n{2}=\lambda\n{1}$ for some real $\lambda$. Hence, in addition to $\n{2}\cdot\p02=0$, one also has $\n{2}\cdot\p01=\lambda\n{1}\cdot\p01=0$. Now Lemma \ref{lem:R.2} implies $\n{2}=\zero$, which is a contradiction.
\end{proof}

\begin{proof}[Proof of Lemma \ref{lem:R.2a}]
Observe that $f:=\ell\cap\conv\P$ is a face of $\conv\P$; see page~162 in \cite{rock}, especially the bottom paragraph there. Hence, by \cite[Corollary 18.1.1]{rock}, $f$ is a closed set. Note next that $f\subseteq\ell$. Also, $f$ is bounded (since $\P$ is so) and convex. It follows that $f=[P,Q]$ for some $P$ and $Q$ in $\conv\P$. Then points $P$ and $Q$ are extreme points of face $f$, and hence of $\conv\P$. Finally, by \cite[Corollary~18.3.1]{rock}, $\{P,Q\}\subseteq\{V_0,\dots,V_{n-1}\}$, so that $\{P,Q\}=\{V_i,V_k\}$ for some $i$ and $k$ in $\intr0{n-1}$.
\end{proof}

\begin{proof}[Proof of Lemma \ref{lem:R.3}]
Assume that, to the contrary,  $V_{j-1}\notin\ext\conv\P$, while conditions (i)--(v) hold. By Lemma \ref{lem:R.2a},
$$f:=\ell\cap\conv\P=[V_i,V_k]$$
for some $i$ and $k$ in $\intr0{n-1}$ such that $V_i$ and $V_k$ are extreme points of $\conv\P$. 

Moreover, $V_{j-1}\notin\ext f$ (since $V_{j-1}$ was assumed to be not an extreme point of $\conv\P$). 
On the other hand, $V_{j-1}\in\ell\cap\conv\P=[V_i,V_k]$.
Therefore,
$$V_{j-1}\in(V_i,V_k)\subseteq\ell.$$
Because $V_j\notin\ell$, one must have $V_{j-1}\ne V_j$. Let then $\ell_1:=V_{j-1}V_j$, the unique line containing the points $V_{j-1}$ and $V_j$. Then, by Lemma \ref{lem:prop:def}, $\ell_1$ is a supporting line to $\P$ and hence to $\conv\P$. It follows that
$f_1:=\ell_1\cap\conv\P$ is a face of $\conv\P$, which contains the point $V_{j-1}$ lying in the relative interior $(V_i,V_k)$ of the segment $[V_i,V_k]$. By the definition of a face, now one has $\{V_i,V_k\}\subset f_1\subseteq\ell_1$. 
Also, the condition $V_{j-1}\in(V_i,V_k)$ implies that $V_i\ne V_k$.
Hence and because $\{V_i,V_k\}\subseteq\ell$, one has $\ell\subseteq\ell_1$, and so, $\ell=\ell_1$. It follows that $V_j\in\ell_1=\ell$, which contradicts the condition $V_j\notin\ell$.   
\end{proof}

\begin{proof}[Proof of Lemma \ref{lem:extreme}]
Suppose that, to the contrary, points $V_0$, $V_1$, and $V_2$ are collinear, while 
$V_1\in\ext\conv\P$. The relation $V_1\in[V_0,V_2]$ is impossible, because $V_1\in\ext\conv\P$ and $\P$ is ordinary. Hence, w.l.o.g., one has 
$$V_0\in(V_1,V_2)$$
(the case $V_2\in(V_1,V_0)$ is quite similar). 
By Lemma \ref{lem:prop:def}, $\ell:=V_1V_2$
is a line supporting to $\P$ and hence also to $\conv\P$. 
Consider the face 
$$f:=\ell\cap\conv\P$$
of $\conv\P$. 
By Lemma \ref{lem:R.2a}, 
$f=[V_i,V_m]$
for some $i,m\in\intr0{n-1}$
such that $V_i$ and $V_m$ are extreme points of $\conv\P$ and hence of face $f$. Therefore, $V_1\in\{V_i,V_m\}$, whence $1\in\{i,m\}$ (since $\P$ is ordinary). Thus, w.l.o.g., $i=1$. Because 
$$V_0\in(V_1,V_2)\subseteq f=[V_1,V_m],$$
it follows that $m\notin\{0,1\}$, and so, 
$$m\in\intr2{n-1}.$$

Consider the set 
$$I:=\{i\in\intr3{n-1}\colon V_i\notin[V_1,V_m]\}.$$
Note that $I$ is non-empty; indeed, otherwise one would have $V_i\in[V_1,V_m]$ for all $i\in\intr0{n-1}$ (because $V_0\in(V_1,V_2)\subseteq[V_1,V_m]$); this would contradict the condition that $\P$ is of dimension 2. 
Hence, 
$$j:=\min I$$
is correctly defined, and then one has $j\in\intr3{n-1}$. Also, $V_{j-1}\in[V_1,V_m]\subseteq\ell$ (if $j=3$ then $V_{j-1}=V_2\in[V_1,V_m]$; if $j\in\intr4{n-1}$ but $V_{j-1}\notin[V_1,V_m]$, then $j-1\in\intr3{n-1}$, whence $j-1\in I$, which contradicts the condition $j=\min I$).
Moreover, $V_j\notin[V_1,V_m]=f=\ell\cap\conv\P$, so that $V_j\notin\ell$ (because $V_j\in\conv\P$). 

Now, by Lemma \ref{lem:R.3}, $V_{j-1}\in\ext\conv\P$. Also, the condition
$j\in\intr3{n-1}$ yields
$j-1\ne1$ and hence $V_{j-1}\ne V_1$ (because $\P$ is ordinary). 
The conditions (i) $V_{j-1}\in[V_1,V_m]$, (ii) $V_{j-1}\ne V_1$, and (iii) $V_{j-1}$ is an extreme point of $\conv\P$ (and hence of face $f=[V_1,V_m]$) imply that $V_{j-1}=V_m$, and so,
$$j-1=m.$$

Consider now the set
$$J:=\{i\in\intr{m+2}n\colon V_i\in\ell\}.$$
Note that $J\ne\emptyset$, because $V_n=V_0\in(V_1,V_2)\subseteq\ell$ and $m+2=j+1\le n$ (since $j\in\intr3{n-1}$\,). 
Hence, the number
$$k:=\min J$$
is correctly defined. Moreover,
$k\in\intr{m+2}n$, $V_k\in\ell$, and $V_{k-1}\notin\ell$ 
(indeed, if $k=m+2$ then $V_{k-1}=V_{m+1}=V_j\notin\ell$; 
if $k\in\intr{m+3}n$ but $V_{k-1}\in\ell$, then $k-1\in\intr{m+2}n$, and so, $k-1\in J$, which contradicts the condition $k=\min J$).
Hence, by Lemma \ref{lem:R.3} (applied to the ``reversed" polygon $(V_{n-1},\dots,V_0)$ in place of $\P=(V_0,\dots,V_{n-1})$), one concludes that $V_k$ is an extreme point---of $\conv\P$ and hence of face $f=[V_1,V_m]$. 
Thus, $V_k\in\{V_1,V_m\}$, and so, 
$$k\in\{1,m\}$$
(since $\P$ is ordinary). 
But $k\in\intr{m+2}n$ and $m=j-1\in\intr2{n-1}$, so that $k\ne m$ and $k\ge4$. This contradicts the conclusion $k\in\{1,m\}$.
\end{proof}

\begin{proof}[Proof of Lemma \ref{lem:case Val extreme}]

Let $\P=(V_0,\dots,V_{n-1})$ be an ordinary quasi-convex polygon. 

By Lemma~\ref{lem:dim=1}, w.l.o.g. $\dim\P=2$, so that $n\ge3$. 

Also, w.l.o.g. $\al=1$, because the property of being an ordinary and quasi-convex polygon is invariant with respect to any cyclic permutation of the indices $(0,\dots,n-1)$.
Then $V_1\in\ext\conv\P$. 

It is enough to show that the sub-polygon $(V_0,V_2,\dots,V_{n-1})$ of $\P$ is quasi-convex, 
because any sub-polygon of an ordinary polygon is obviously ordinary.

By Lemma \ref{lem:prop:def}, it is enough to prove that $V_3,\dots,V_{n-1}\|[V_0,V_2]$. 

By Lemma \ref{lem:extreme},
the vertices $V_0$, $V_1$, and $V_2$ are non-collinear. 

Therefore, w.l.o.g., 
$$V_1=(0,0),\quad V_0=(1,0),\quad V_2=(0,1).$$
Indeed, on the one hand, the notion of quasi-convexity is invariant under one-to-one affine transformations of $\R^2$ onto itself, and, on the other hand, by \cite[Theorem 1.6]{rock}, there exists a one-to-one affine transformation of $\R^2$ which carries any three given non-collinear points
to points $(0,0)$, $(1,0)$, and $(0,1)$. 

Introduce the vector
$$\nn:=(1,1).$$
Then $\nn\ne\zero$ and
$\nn\cdot\v02=0$.
Therefore, to verify the relation $V_3,\dots,V_{n-1}\|[V_0,V_2]$, it suffices to show that 
$$\nn\cdot\v0\beta \ge0,\quad
\text{for every } 
\beta \in\intr3{n-1};$$ 
moreover, in view of \cite[Corollary 18.5.1]{rock}, one may assume w.l.o.g. that $V_\beta\in\ext\conv\{V_0,V_2,\dots,V_{n-1}\}$. 
Then, by Lemma \ref{lem:extreme}, the points 
$$V_{\beta-1},V_\beta, V_{\beta+1}\quad\text{are non-collinear}.$$ 

Now let $(x,y)$ denote the coordinates of such a point $V_\beta\in\ext\conv\{V_0,V_2,\dots,V_{n-1}\}$:
$$V_\beta =(x,y).$$
Since $\P$ is quasi-convex, by Lemma \ref{lem:prop:def} one has $V_2,V_\beta \|[V_1,V_0]$. Because the line $V_1V_0$ is the set $\{(u,v)\in\R^2\colon v=0\}$, the corresponding closed half-plane containing point $V_2=(0,1)$ must be $\{(u,v)\in\R^2\colon v\ge0\}$, which implies that 
$$y\ge0,$$
because the half-plane must contain the point $V_\beta =(x,y)$ as well. Similarly,
$$x\ge0.$$ 

Since $\v\beta0 =(1-x,-y)$, $\v\beta1 =(-x,-y)$, $\v\beta2 =(-x,1-y)$, and $\nn\cdot\v0\beta =x+y-1$, it is straightforward to verify the identity
\begin{equation}\label{eq:identity}
x\,\v\beta0 +y\,\v\beta2 =(\nn\cdot\v0\beta )\,\v\beta1 .
\end{equation} 

Because $\P$ is quasi-convex, Lemma \ref{lem:prop:def} implies that there exist nonzero vectors $\n{\beta} $ and $\m{\beta} $ such that 
\begin{gather}
\n{\beta} \cdot\v \beta {\beta -1}=0,\quad\m{\beta} \cdot\v \beta {\beta +1}=0,\quad \text{while} \notag\\
\n{\beta} \cdot\v \beta j\ge0,\quad\m{\beta} \cdot\v \beta j\ge0\quad
\text{for all $j\in\intr0{n-1}$. }\label{eq:ge0}
\end{gather}

Assume that $\nn\cdot\v0\beta<0$, which is the contrary to what we must prove. Dot-multiply both sides of identity \eqref{eq:identity} by $\n{\beta}$ and $\m{\beta}$. Then the inequalities $x\ge0$, $y\ge0$, \eqref{eq:ge0}, and $\nn\cdot\v0\beta<0$ imply that 
$$\n{\beta}\cdot\v\beta1=0\quad\text{and}\quad\m{\beta}\cdot\v\beta1=0.$$
Now Lemma \ref{lem:R.1} (with $V_\beta$, $V_{\beta-1}$, $V_{\beta+1}$, $\n{\beta}$, $\m{\beta}$, $\v\beta1$ in place of, respectively, $P_0$, $P_1$, $P_2$, $\n1$, $\n2$, $\aa$) 
yields $\v\beta1=\zero$, that is, $V_\beta=V_1$, which contradicts the conditions that $\beta \in\intr3{n-1}$ and $\P$ is ordinary. 
\end{proof}

\begin{proof}[Proof of Lemma \ref{lem:R.13}]
Assume the contrary: that 
$$Q:=\{q\in\intr1{k-1}\colon V_q\notin\ell\}\ne\emptyset\quad\text{and}\quad
R:=\{q\in\intr1{k-1}\colon V_q\in\ell\}\ne\emptyset.$$
Then
$$\tau:=\min Q\quad\text{and}\quad\si:=\min R$$
are correctly defined. 
At that $\tau\in\intr1{k-1}$ and $V_\tau\notin\ell$. Also, $V_{\tau-1}\in\ell$. (Indeed, if $\tau=1$ then $V_{\tau-1}=V_0\in[V_0,V_k]=\ell\cap\conv\P\subseteq\ell$; 
and if $\tau\ge2$ then $\tau-1\in\intr1{k-1}$, so that $V_{\tau-1}\notin\ell$ would imply $\tau-1\in Q$, which would contradict the condition $\tau=\min Q$.)

Now, by Lemma \ref{lem:R.3}, $V_{\tau-1}$ is an extreme point of $\conv\P$ and hence of face
$\ell\cap\conv\P=[V_0,V_k]$. 
Therefore, $\tau-1\in\{0,k\}$ (since $\P$ is ordinary). But $\tau-1\ne k$, because $\tau\in\intr1{k-1}$. It follows that $\tau-1=0$, so that $\tau=1$, whence $V_1=V_\tau\notin\ell$. This means that $1\notin R$, so that 
$\si\in\intr2{k-1}$. Also, $V_\si\in\ell$.
Moreover, $V_{\si-1}\notin\ell$. (Indeed, the condition $\si\in\intr2{k-1}$ yields $\si-1\in\intr1{k-1}$, so that $V_{\si-1}\in\ell$ would imply $\si-1\in R$, which would contradict the condition $\si=\min R$.)

By Lemma \ref{lem:R.3} (applied to the ``reversed" polygon $(V_{n-1},\dots,V_0)$), $V_\si$ is an extreme point of $\conv\P$ and hence of $\ell\cap\conv\P=[V_0,V_k]$. Therefore $\si\in\{0,k\}$, which contradicts the condition $\si\in\intr2{k-1}$.  
\end{proof}

\begin{proof}[Proof of Lemma \ref{lem:prop:conv-reduction}]
Suppose that $V_\al$ is not an extreme point of the convex hull of an ordinary polygon $\P=(V_0,\dots,V_{n-1})$, for some $\al\in\intr0{n-1}$. 
Then the claim $\conv\P^{(\al)}=\conv\P$ follows immediately from the definitions of the convex hull and an extreme point.

It remains to prove that $\edg\P^{(\al)}=\edg\P$.
By Lemma \ref{lem:prop:def}, the line $\ell:=V_\al V_{\al+1}$ is supporting to $\conv\P$. In view of Lemma \ref{lem:R.2a}, w.l.o.g. $\ell\cap\conv\P=[V_0,V_k]$ for some $k\in\intr1{n-1}$, and $V_0$ and $V_k$ are extreme points of $\conv\P$.

Since $V_\al\notin\ext\conv\P$, one has $V_\al\notin\{V_0,V_k\}$, whence $\al\notin\{0,k\}$ (since $\P$ is ordinary). That is,
$$\al\in\intr1{k-1}\cup\intr{k+1}{n-1}.$$

Thus, one has only these two cases: $\al\in\intr1{k-1}$ and $\al\in\intr{k+1}{n-1}$. These two cases are quite similar. In fact, say, the latter case can be reduced to the former one by considering the cyclic permutation $(V_k,\dots,V_{n-1},V_0,\dots,V_{k-1})$ of the vertices of polygon $\P=(V_0,\dots,V_{n-1})$, which preserves the union of the edges. 
Thus, w.l.o.g.,
$$\al\in\intr1{k-1}.$$

Now it follows that $\forall q\in\intr1{k-1}\ V_q\in\ell$. 
Therefore, $\forall q\in\intr1{k-1}\ V_q\in\ell\cap\conv\P=[V_0,V_k]$. 
Hence, by Lemma \ref{lem:R.5}, 
$$\bigcup_{q\in\intr0{k-1}}[V_q,V_{q+1}]=[V_0,V_k];$$
also, using again Lemma \ref{lem:R.5} and taking into account the condition $\al\in\intr1{k-1}$, one has 
$$\bigcup_{q\in\intr0{\al-2}}[V_q,V_{q+1}]
\cup[V_{\al-1},V_{\al+1}]\cup\bigcup_{q\in\intr{\al+1}{k-1}}[V_q,V_{q+1}]
=[V_0,V_k],$$
whence
\begin{align*}
\edg\P^{(\al)}& \\
=&\bigcup_{q\in\intr0{\al-2}}[V_q,V_{q+1}]
\cup[V_{\al-1},V_{\al+1}]\cup\bigcup_{q\in\intr{\al+1}{k-1}}[V_q,V_{q+1}] 
\cup\bigcup_{q\in\intr{k}{n-1}}[V_q,V_{q+1}] \\
=&[V_0,V_k]\cup\bigcup_{q\in\intr{k}{n-1}}[V_q,V_{q+1}] \\
=&\bigcup_{q\in\intr0{k-1}}[V_q,V_{q+1}]\cup\bigcup_{q\in\intr{k}{n-1}}[V_q,V_{q+1}]
=\edg\P.
\end{align*}
\end{proof}

\begin{proof}[Proof of Lemma \ref{lem:conv-reduction}]
This follows immediately from Lemmas \ref{lem:case Val extreme} and \ref{lem:prop:conv-reduction}.
\end{proof}

\begin{proof}[Proof of Lemma~\ref{lem:new}]
Suppose, to the contrary, that $n\ge5$ and $\P=(V_0,\dots,V_{n-1})$ is a quasi-convex polygon such that $\P^{(\al)}$ is convex for every $\al\in\intr0{n-1}$, while polygon $\P$ is not convex. 
Then there is a point 
$$
P\in(\partial\conv\P)\setminus\edg\P.
$$
By \cite[Theorem 11.6]{rock}, there is a line $\ell$ containing $P$ and supporting to $\conv\P$. 
By Lemma~\ref{lem:R.2a}, w.l.o.g.\ 
$$\ell\cap\conv\P=[V_0,V_k]\ni P$$
for some $k\in\intr0{n-1}$, and $V_0$ and $V_k$ are extreme points of $\conv\P$.
In particular, this implies that $V_k\ne V_1$ and hence
$k\ne1$ (because otherwise one would have $P\in[V_0,V_1]$ while
$P\notin\edg\P$). 

Also, $V_{i+1}\ne V_i$ for any $i\in\intr0{n-1}$. Indeed, otherwise $\conv\P^{(i+1)}=\conv\P$ and $\edg\P^{(i+1)}=\edg\P$, which contradicts the assumptions that $\P^{(\al)}$ is convex for every $\al\in\intr0{n-1}$ and $\P$ is not convex. 

In particular, $V_1\ne V_0$, so that $V_1\notin\{V_0,V_k\}$.

Hence, in the case when $V_1\in\ell$ one has $V_1\in\ell\cap\conv\P=[V_0,V_k]$, and so, $V_1\in(V_0,V_k)$ and hence $V_1\notin\ext\conv\P$.
In view of Lemma~\ref{lem:R.3} (with $j=2$), it follows that 
$V_2\in\ell$ and hence
$V_2\in\ell\cap\conv\P=[V_0,V_k]$.
Therefore, $\conv\P^{(1)}=\conv\P$, whence $\partial\conv\P^{(1)}=\partial\conv\P$, while
$[V_0,V_2]\subseteq[V_0,V_1]\cup[V_1,V_2]$, so that $$\edg\P^{(1)}\subseteq\edg\P\subsetneq\partial\conv\P
=\partial\conv\P^{(1)},$$
which contradicts the condition that $\P^{(\al)}$ is convex for every $\al\in\intr0{n-1}$.

It remains to consider the case when $V_1\notin\ell$. 
One has $P\in[V_0,V_k]\subseteq\conv\P^{(1)}$ (since $k\ne1$), $P\in\partial\conv\P$, and $\interi\conv\P^{(1)}\subseteq\interi\conv\P$. Hence,
$P\in\partial\conv\P^{(1)}$, so that $P\in\edg\P^{(1)}$ (since $\P^{(1)}$ is convex), while $P\notin\edg\P$. So,
$$P\in(V_0,V_2).$$
In particular, this implies that $P$ is not an extreme point of $\conv\P$.
Since $P\in[V_0,V_k]$ and $V_0$ and $V_k$ are extreme points of $\conv\P$, one has 
$P\in(V_0,V_k)\subseteq\ell$. 
It follows that 
$$\ell=V_0V_2=V_0V_k,$$
and so, $\P$ is to one side of $[V_0,V_2]$. 

Since $\P$ is quasi-convex, it is also to one side of $[V_0,V_1]$ and to one side of $[V_1,V_2]$ (by Lemma~\ref{lem:prop:def}).
The condition $V_1\notin\ell$ and the conclusion that $\ell=V_0V_2$ imply that $\dim\conv\{V_0,V_1,V_2\}=2$. By \cite[Theorem 1, page 31]{grun}, $\conv\{V_0,V_1,V_2\}$ is the intersection of three closed half-planes whose boundaries are the lines $V_0V_1$, $V_1V_2$, and $V_0V_2$. Since $\P=(V_0,\dots,V_{n-1})$ is to one side of each of the segments $[V_0,V_1]$, $[V_1,V_2]$, and $[V_0,V_2]$, it follows that $\conv\P\subseteq\conv\{V_0,V_1,V_2\}$ and hence $\conv\P=\conv\{V_0,V_1,V_2\}$.

Moreover, $\{V_0,\dots,V_{n-1}\}\subseteq\edg\P\subseteq\partial\conv\P$, since $\P$ is quasi-convex. Thus,
$\{V_0,\dots,V_{n-1}\}\subseteq\partial\conv\{V_0,V_1,V_2\}=[V_0,V_1]\cup[V_1,V_2]\cup[V_0,V_2]$.

Furthermore, $\{V_0,\dots,V_{n-1}\}\subseteq\{V_0,V_1,V_2\}.$ Indeed, otherwise there is some $\al\in\intr0{n-1}$ such that 
$V_\al\in(V_0,V_1)\cup(V_1,V_2)\cup(V_0,V_2)$; then in fact necessarily $\al\in\intr3{n-1}$, and w.l.o.g.\ 
$V_\al\in(V_0,V_1)$ \big(the cases $V_\al\in(V_1,V_2)$ and $V_\al\in(V_0,V_2)$ are quite similar\big).
By Lemma~\ref{lem:R.3} (with $j=\al+1$), one has 
$V_{\al+1}\in V_0V_1$ (since the line $V_0V_1$ is supporting to $\conv\P$);
similarly, applying the same Lemma~\ref{lem:R.3} to the ``reversed'' polygon $(V_0,V_{n-1},\dots,V_1)$, one sees that 
$V_{\al-1}\in V_0V_1$.
Therefore, $[V_{\al-1},V_{\al+1}]\subseteq[V_{\al-1},V_\al]\cup[V_\al,V_{\al+1}]$. Hence and because $\al\ge3$,  $$\edg\P^{(\al)}\subseteq\edg\P\subsetneq\partial\conv\P=\partial\conv\{V_0,V_1,V_2\}\subseteq\partial\conv\P^{(\al)},$$
which contradicts the condition that $\P^{(\al)}$ is convex for every $\al\in\intr0{n-1}$. 

Next, $\{V_i,V_{i+1}\}\ne\{V_0,V_2\}$ for any $i\in\intr0{n-1}$. Indeed, otherwise $P\in[V_0,V_2]=[V_i,V_{i+1}]\subseteq\edg\P$, which contradicts assumption $P\in(\partial\conv\P)\setminus\edg\P$. 

Recall also that $V_{i+1}\ne V_i$ for any $i\in\intr0{n-1}$, as was shown above in this proof.

Thus, polygon $\P=(V_0,V_1,V_2,\dots,V_{n-1})$ is a sequence of vertices of length $n\ge5$ such that (i) every vertex of $\P$ equals $V_0$, $V_1$, or $V_2$; (ii) (every vertex which equals) $V_0$ is followed by $V_1$; (iii) $V_1$ is followed by $V_0$ or $V_2$; (iv) $V_2$ is followed by $V_1$; and (v) $V_{n-1}=V_1$ (since $V_n=V_0$).

It follows that 
$(V_0,V_1,V_2,V_3,V_4,V_5)$ equals either $(V_0,V_1,V_2,V_1,V_0,V_1)$ or \break
$(V_0,V_1,V_2,V_1,V_2,V_1)$ (and necessarily $n\ge6$). 
Then $\P^{(4)}=(V_0,V_1,V_2,V_1,V_1,\dots,\\
V_{n-1})$, so that $\conv\P^{(4)}=\conv\P$ and $\edg\P^{(4)}=\edg\P$, which contradicts the assumption that $\P^{(\al)}$ is convex for every $\al\in\intr0{n-1}$ while $\P$ is not convex. 
\end{proof}

\begin{proof}[Proof of Lemma \ref{lem:R.10}]
Suppose the contrary, that the conditions listed in Lemma~\ref{lem:R.10} can all be satisfied at once.
Since $\P$ is ordinary, all the $V_i$'s must be distinct. Since $\{V_0,V_2\}\subseteq\ell$ and $\{V_1,V_3\}\cap\ell=\emptyset$, the points $V_0$, $V_1$, and $V_2$ must be non-collinear. Hence,
w.l.o.g. $V_0=(0,0)$, $V_1=(0,1)$, and $V_2=(1,0)$. 
Let $V_3=(x,y)$, for some real $x$ and $y$. 

Introduce vectors $\n1:=(1,0)$ and $\n2:=(0,1)$. 
Then $V_1,V_3\|[V_0,V_2]$, since line $\ell$ is supporting to $\conv\P$ and contains points $V_0$ and $V_2$. 
Also, $\n2\cdot\v02=0<1=\n2\cdot\v01$. 
Hence,
$y=\n2\cdot\v03\ge0$. Moreover, 
$$y>0,$$
because $V_3\notin\ell$. 

Next, $\n1\cdot\v01=0<1=\n1\cdot\v02$. Also, by Lemma \ref{lem:prop:def}, $V_2,V_3\|[V_0,V_1]$. Hence, 
$$x=\n1\cdot\v03\ge0.$$

Similarly, using the vector $\nn:=(-1,-1)$ and condition $V_0,V_3\|[V_1,V_2]$, one has 
$\nn\cdot\v12=0<1=\nn\cdot\v10$, whence 
$$1-x-y=\nn\cdot\v13\ge0.$$

Further, using the vector $\mm:=(-y,x-1)$ and condition $V_0,V_1\|[V_2,V_3]$, one has 
$\mm\cdot\v23=0<y=\mm\cdot\v20$, whence 
$x+y-1=\mm\cdot\v21\ge0$.  
Thus, $x+y-1\ge0$ while $1-x-y\ge0$; that is, 
$$x+y=1.$$

Finally, using the vector $\ee:=(y,-x)$, one has 
$\ee\cdot\v30=0$, $\ee\cdot\v31=-x\le0$, and $\ee\cdot\v32=y>0$. Since  $V_1,V_2\|[V_0,V_3]$, one must have $-x\ge0$, so that $x=0$. Now $x+y=1$ yields $y=1$, so that $V_3=(x,y)=(0,1)=V_1$, which contradicts the condition that $\P$ is ordinary.  
\end{proof}

\begin{proof}[Proof of Lemma \ref{lem:R.11}]
Suppose the contrary, that the conditions listed in Lemma \ref{lem:R.11} can all be satisfied at once.
Then condition $k\in\intr2{n-2}$ implies $n\ge4$.
If $n=4$, then $k\in\intr2{n-2}$ also implies $k=2$, and so, in this case Lemma \ref{lem:R.11} follows from Lemma \ref{lem:R.10}.

The remaining case, $n\ge5$, can be proved by induction in $n$. Indeed, if $n\ge5$, then at least one of the two sets, $\{V_1,\dots,V_{k-1}\}$ or $\{V_{k+1},\dots,V_{n-1}\}$, contains at least two distinct vertices. W.l.o.g., the first one of these two quite similar subcases takes place; that is, $k-1\ge2$, so that 
$$k\in\intr3{n-2}.$$ 

Consider now the polygon $\P^{(1)}=(V_0,V_2,V_3,\dots,V_{n-1})$. 
Then
$$\P^{(1)}=(U_0,\dots,U_{n-2}),$$
where
$$U_q:=
\begin{cases}
V_0 & \text{ if }\ q=0, \\
V_{q+1} & \text{ if }\ q\in\intr1{n-2}.
\end{cases}
$$
By Lemma \ref{lem:conv-reduction}, $\P^{(1)}$ is an ordinary quasi-convex polygon.
The line $\ell$ is supporting to $\conv\P^{(1)}$, because $\ell$ is supporting to $\conv\P$ and $V_0\in\ell\cap\P^{(1)}$.
Next, 
$$\ell\cap\conv\P^{(1)}\supseteq\ell\cap\P^{(1)}\supseteq\{V_0,V_k\},$$
because $\ell\supseteq\{V_0,V_k\}$ and $k\ge3$. 
Hence,
$$\ell\cap\conv\P\supseteq\ell\cap\conv\P^{(1)}\supseteq[V_0,V_k]=\ell\cap\conv\P,$$
which yields 
$$\ell\cap\conv\P^{(1)}=[V_0,V_k]=[U_0,U_{\tilde k}],$$
where
$$\tilde k:=k-1\in\intr2{(n-1)-2}$$
(because $k\in\intr3{n-2}$). 
Finally, 
$$\{U_1,\dots,U_{\tilde k-1}\}\cap\ell=\{V_2,\dots,V_{k-1}\}\cap\ell
\subseteq\{V_1,\dots,V_{k-1}\}\cap\ell=\emptyset$$
and
$\{U_{\tilde k+1},\dots,U_{(n-1)-1}\}\cap\ell=\{V_{k+1},\dots,V_{n-1}\}\cap\ell
=\emptyset$.
Hence, $\P^{(1)}$, $\tilde k$, and $n-1$ satisfy (in place of $\P$, $k$, and $n$, respectively) all the conditions listed in Lemma \ref{lem:R.11}. 
Thus, the induction step is verified. 
\end{proof}

\begin{proof}[Proof of Lemma \ref{lem:conv1}]
If a polygon is convex, then it is trivially quasi-convex.

Assume now that, vice versa, a polygon $\P:=(V_0,\dots,V_{n-1})$ is ordinary and quasi-convex. Take any point $V\in\partial\conv\P$. What we have then to show is that $V\in[V_q,V_{q+1}]$ for some $q\in\intr0{n-1}$. 

Since $V\in\partial\conv\P$, by \cite[Corollary 11.6.1]{rock} there exists a line $\ell$ containing point $V$ and supporting to $\conv\P$, so that $V\in\ell\cap\conv\P$. Then, by Lemma \ref{lem:R.2a}, one has w.l.o.g. that
$$\ell\cap\conv\P=[V_0,V_k]$$
for some $k\in\intr1{n-1}$. 
If $k=1$, then $V\in\ell\cap\conv\P=[V_0,V_k]=[V_q,V_{q+1}]$ for $q=0$.
Similarly, if $k=n-1$ then $V\in\ell\cap\conv\P=[V_0,V_k]=[V_k,V_0]=[V_q,V_{q+1}]$ for $q=n-1$.
It remains to consider the case $k\in\intr2{n-2}$. 
Then, by Lemma \ref{lem:R.11}, either $\{V_1,\dots,V_{k-1}\}\cap\ell\ne\emptyset$ or 
$\{V_{k+1},\dots,V_{n-1}\}\cap\ell\ne\emptyset$. 
These two cases are quite similar to each other.
Hence, w.l.o.g., one has $\{V_1,\dots,V_{k-1}\}\cap\ell\ne\emptyset$, that is, $\exists q\in\intr1{k-1}$ $V_q\in\ell$. 
Now (the second one of the two equivalent implications in) Lemma \ref{lem:R.13} yields $\forall q\in\intr1{k-1}$ $V_q\in\ell$. 
Therefore, in view of Lemma \ref{lem:R.5}, 
$$V\in[V_0,V_k]=\bigcup_{q\in\intr0{k-1}}[V_q,V_{q+1}],$$
so that indeed $V\in[V_q,V_{q+1}]$ for some $q\in\intr0{n-1}$. 
\end{proof}

\begin{proof}[Proof of Lemma \ref{lem:th:conv-reduction}]
This follows immediately from Lemma \ref{lem:conv1}
and Lemma \ref{lem:conv-reduction}.
\end{proof}

\begin{proof}[Proof of Lemma \ref{lem:cor:sub-polygon}]
This follows immediately from Lemma \ref{lem:th:conv-reduction}.
\end{proof}

\begin{proof}[Proof of Lemma \ref{lem:sep}]
This follows immediately from \cite[Theorem 1, page 31]{grun}.
\end{proof}

\begin{proof}[Proof of Lemma \ref{lem:Q}]
Let $H$ be a closed half-plane such that $\partial H=\ell$ and $H\supseteq\conv\{U_0,U_1,U_2\}$. Then 
$H\supseteq\conv\{U_0,U_2\}=[U_0,U_2]\ni O$. 
On the other hand,
(i) $O\in(U_1,U_3)$ implies that $\overrightarrow{U_1U_3}=\la\overrightarrow{U_1O}$ for some $\la>0$ (in fact, $\la>1$) and 
(ii) $U_1\in\ell=\partial H$ implies that the half-plane $H$ is a cone with vertex $U_1$. 
Hence, $U_3\in H$, so that $\conv\Q=\conv\{U_0,U_1,U_2,U_3\}\subseteq H$.
\end{proof}

\begin{proof}[Proof of Lemma~\ref{lem:prop:apart}]
Suppose that the conditions of Proposition \ref{prop:apart} hold. 
Then $\Q:=(U_0,U_1,U_2,U_3):=(V_0,V_\al,V_i,V_\beta)$ is a sub-$4$-gon of $\P$ and hence strictly convex, in view of Lemma~\ref{lem:cor:sub-polygon}. 
Let $\ell:=U_0U_2$. If $V_\al,V_\beta\|[V_0,V_i]$ (i.e., $U_1,U_3\|[U_0,U_2]$), then the line $\ell$ is supporting to $\conv\Q$. 
It follows from Lemma~\ref{lem:R.2a} and the strictness of $\Q$ that $\ell\cap\conv\Q=[U_0,U_2]$. 
Now Proposition~\ref{prop:apart} follows from Lemma~\ref{lem:R.11} (with $n=4$ and $k=2$) and the strictness of $\Q$.
\end{proof}

\begin{proof}[Proof of Lemma \ref{lem:if}]
This follows immediately from the definitions.
\end{proof}

\begin{proof}[Proof of Lemma \ref{lem:0213}]
By Lemma~\ref{lem:prop:apart}, the strict convexity of $\Q$ implies \break
$U_1U_3\notOneSide[U_0U_2]$. It remains to refer to Lemma~\ref{lem:prop:def}.
\end{proof}

\begin{proof}[Proof of Lemma \ref{lem:si(1)}]
Suppose that, on the contrary, $1\si\ne1$, while all the conditions of Lemma~\ref{lem:si(1)} hold. Then $1\si\in\intr2{n-1}$, because $1\si\ne0\si=0$. Let $i:=1\si^{-1}$, so that $i\si=1$. Then $i\ne0$, since $0\si=0$. Also, $i\ne1$, because of the assumption $1\si\ne1$. Finally,
$i\ne n-1$, because $i\si=1<1\si<(n-1)\si$. Hence, $i\in\intr2{n-2}$.
It also follows that
$0\si<i\si<1\si<(n-1)\si$.
Therefore, $\Q:=(V_{0\si},V_{i\si},V_{1\si},V_{(n-1)\si})$ is a sub-$4$-gon of the strictly convex polygon $\P$, and so, 
$\Q$ is strictly convex,
in view of Lemma~\ref{lem:cor:sub-polygon}.  

On the other hand, $\tilde\Q:=(V_{0\si},V_{1\si},V_{i\si},V_{(n-1)\si})$ is a sub-$4$-gon of the strictly convex polygon $\P\si$ (because $0<1<i<n-1$), so that $\tilde\Q$ is also strictly convex. This contradicts Lemma \ref{lem:0213}.
\end{proof}

\begin{proof}[Proof of Lemma \ref{lem:si}]
Suppose that all the conditions of Lemma~\ref{lem:si} hold. 
Let
$$J:=\{j\in\intr0{n-1}\colon i\si=i\ \forall i\in\intr0j\}.$$
By Lemma \ref{lem:si(1)}, $\{0,1\}\subseteq J$. Let
$$k:=\max J,$$
so that $k\in\intr1{n-1}$. 
It suffices to show that $k=n-1$. 
Assume the contrary: $k\in\intr1{n-2}$. Let
$$\pi:=\th^k\si\th^{-k} \quad\text{and}\quad
\Q:=\P\th^{-k},$$
where $\th$ is the ``primitive" cyclic permutation defined before the statement of Proposition \ref{prop:G}. 
By Lemma \ref{lem:if}, polygons $\Q=\P\th^{-k}$ and $\Q\pi=(\P\si)\th^{-k}$ are strictly convex, because polygons $\P$ and $\P\si$ are strictly convex and $\th^{-k}\in G_n$. Also, recalling the definition $k=\max J$, one has
$$
\begin{aligned}
0\pi&=0\th^k\si\th^{-k} =k\si\th^{-k}=k\th^{-k}=0,\\  
(n-1)\pi&=(n-1)\th^k\si\th^{-k} =(k-1)\si\th^{-k}=(k-1)\th^{-k}=n-1.
\end{aligned}
$$
Next, $k\in\intr1{n-2}$ implies $n\ge3$, so that $1\ne n-1$ and hence $1\pi\ne(n-1)\pi$. 
Therefore, $1\pi\in\intr0{n-1}\setminus\{(n-1)\pi\}=\intr0{n-1}\setminus\{n-1\}
=\intr0{n-2}$, so that $1\pi<n-1=(n-1)\pi$.
Applying now Lemma \ref{lem:si(1)} (with $\Q$ and $\pi$ in place of $\P$ and $\si$), one has
$1=1\pi=1\th^k\si\th^{-k}$, whence $1\th^k=1\th^k\si$.
Since $1\th^k=k+1$, one now has 
$(k+1)\si=k+1$. 
Thus, $k+1\in J$, which contradicts the condition $k=\max J$. 
\end{proof}

\begin{proof}[Proof of Lemma~\ref{lem:al-beta}]
W.l.o.g., $0=\al<\beta\le n-1$. Note that $n\ge3$ and $\P$ is ordinary, since $\P$ is strict. Let $\ell_0:=V_\al V_{\al+1}=V_0V_1$. Then, by Lemma~\ref{lem:prop:def}, the line $\ell_0$ is supporting to $\conv\P$, whence $\{V_\beta,V_{\beta+1}\}\subseteq H$ for a closed half-plane $H$ with $\partial H=\ell_0$.

If $\{V_\beta,V_{\beta+1}\}\subseteq\ell_0=V_0V_1$, then $\beta\in\{0,1\}$ and $\beta+1\in\{1,n\}$ (since $\P$ is strict). 
But $\beta\in\intr1{n-1}$ and hence $\beta\ne0$. Now if $\beta=1$, then $\beta+1\notin\{1,n\}$, since $n\ge3$. This contradiction shows that $\{V_\beta,V_{\beta+1}\}\not\subseteq\ell_0$.

Therefore, one has one of the following three cases.

\emph{Case 1}: $\{V_\beta,V_{\beta+1}\}\cap\ell_0=\emptyset$.\quad
Then $\{V_\beta,V_{\beta+1}\}\subseteq H\setminus\ell_0=\interi H$, whence $[V_\beta,V_{\beta+1})\subseteq[V_\beta,V_{\beta+1}]\subseteq\interi H$, while $[V_0,V_1)\subseteq\ell_0=\partial H$, which implies the conclusion of Lemma~\ref{lem:al-beta} (with $\al=0$). 

\emph{Case 2}: $V_\beta\in\ell_0$ and $V_{\beta+1}\notin\ell_0$.\quad
Then $V_\beta=V_1$ (since $\P$ is strict, $\ell_0=V_0V_1$, and $\beta\in\intr1{n-1}$). 
Let now $\ell_\beta:=V_\beta V_{\beta+1}$. Then $\ell_\beta\ne\ell_0$ (since $V_{\beta+1}\in\ell_\beta\setminus\ell_0$), so that 
$\ell_0\cap\ell_\beta=\{V_\beta\}$. It follows that 
$[V_0,V_1)\cap[V_\beta,V_{\beta+1})\subseteq\ell_0\cap\ell_\beta=\{V_\beta\}$. 
But $V_\beta=V_1\notin[V_0,V_1)$, which implies the conclusion of Lemma~\ref{lem:al-beta} (with $\al=0$) in Case~2 as well.  
 
\emph{Case 3}: $V_\beta\notin\ell_0$ and $V_{\beta+1}\in\ell_0$.\quad
This case is quite similar to Case~2.
Indeed, 
here $V_{\beta+1}=V_0$ (since $\P$ is strict, $\ell_0=V_0V_1$, and $\beta+1\in\intr2n$);
at that, $\beta+1=n$ and hence $\beta=n-1$. 
Consider again $\ell_\beta=V_\beta V_{\beta+1}$. Then $\ell_\beta\ne\ell_0$ (since $V_\beta\in\ell_\beta\setminus\ell_0$), so that 
$\ell_0\cap\ell_\beta=\{V_{\beta+1}\}$. It follows that 
$[V_0,V_1)\cap[V_\beta,V_{\beta+1})\subseteq\ell_0\cap\ell_\beta=\{V_{\beta+1}\}=\{V_0\}$. 
But $V_0\notin[V_{n-1},V_0)=[V_\beta,V_{\beta+1})$, which implies the conclusion of Lemma~\ref{lem:al-beta} (with $\al=0$) in Case~3, too.  
\end{proof}

\bibliographystyle{amsplain}

\end{document}